\documentclass[11pt]{amsart}
\usepackage{a4wide}
\usepackage{amsmath, amsfonts, amssymb,amscd,graphicx}
\usepackage{mathrsfs,url}

\date{\today}

\author[M. Bodirsky]
{Manuel Bodirsky}
    \address{Institut f\"{u}r Algebra\\TU Dresden\\01062 Dresden\\Germany}
    \email{Manuel.Bodirsky@tu-dresden.de}
   \urladdr{http://www.math.tu-dresden.de/~bodirsky/}

\author[M.~Pinsker]
{Michael Pinsker}
	\address{Department of Algebra, MFF UK, Sokolovska 83, 186 00 Praha 8, Czech Republic}        
	\email{marula@gmx.at}
    \urladdr{http://dmg.tuwien.ac.at/pinsker/}

\author[A.~Pongr\'{a}cz]
{Andr\'{a}s Pongr\'{a}cz}
    \address{Department of Algebra and Number Theory\\
    University of Debrecen\\
    4032 Debrecen, Egyetem square 1\\
    Hungary}
\email{pongracz.andras@science.unideb.hu}

\thanks{The first and third author have received funding from the European Research Council under the European Community's Seventh Framework Programme (FP7/2007-2013 Grant Agreement no. 257039). The first author also received funding from the German Science Foundation (DFG, project number 622397). The second author has been funded through projects I836-N23 and P27600 of the  Austrian Science Fund (FWF). The third author was also partially supported by the Hungarian Scientific Research Fund (OTKA) grant no.~K109185.}

\title{Reconstructing the topology of clones}

\DeclareMathOperator{\HHH}{H}
\DeclareMathOperator{\PPP}{P}
\DeclareMathOperator{\PPPfin}{P^{\fin}}

\DeclareMathOperator{\Ima}{Im}

\newcommand{\fC}{\mathscr C}
\newcommand{\fD}{\mathscr D}

\DeclareMathOperator{\fin}{fin}

\DeclareMathOperator{\pol}{Pol}
\DeclareMathOperator{\Emb}{Emb}

\DeclareMathOperator{\Clo}{Clo}

\newcommand{\HSPfin}{\HSP^{\fin}}
\newcommand{\rest}{{\upharpoonright}}

\newcommand{\Fresse}{Fra\"{i}ss\'{e}}

\newcommand{\ignore}[1]{}

\newcommand{\To}{\rightarrow}
\newcommand{\nin}{\notin}

\newcommand{\mult}{\times}

\DeclareMathOperator{\Csp}{CSP}

\DeclareMathOperator{\id}{id}

\DeclareMathOperator{\Aut}{Aut}
\DeclareMathOperator{\End}{End}
\DeclareMathOperator{\Pol}{Pol}

\newcommand{\G}{{\bf G}}

\newcommand{\F}{\mathscr F}
\renewcommand{\H}{\mathscr H}

\newcommand{\M}{\mathscr M}

\theoremstyle{plain}

    \newtheorem{thm}{Theorem}
    \newtheorem{theorem}[thm]{Theorem}

    \newtheorem{lem}[thm]{Lemma}
    \newtheorem{lemma}[thm]{Lemma}

    \newtheorem{prop}[thm]{Proposition}
    \newtheorem{proposition}[thm]{Proposition}

    \newtheorem{cor}[thm]{Corollary}
    \newtheorem{corollary}[thm]{Corollary}

    \newtheorem{quest}{Question}

\theoremstyle{definition}

    \newtheorem{defn}[thm]{Definition}
    \newtheorem{definition}[thm]{Definition}

\DeclareMathOperator{\HSP}{HSP}

\DeclareMathOperator{\SSS}{S}

\DeclareMathOperator{\comp}{comp}

\DeclareMathOperator{\Dom}{Dom}

\DeclareMathOperator{\Age}{Age}

\newtheorem{remark}[thm]{Remark}

\begin{document}
\begin{abstract}
Function clones are sets of functions on a fixed domain 
that are closed under composition and 
contain the projections.
They carry a natural algebraic structure, provided by the laws of composition which hold in them, as well as a natural topological structure, provided by the topology of pointwise convergence, under which composition 
of functions becomes continuous. Inspired by recent results indicating the importance of the topological ego of function clones even for originally algebraic problems,  we study questions of the following type: In which situations does the algebraic structure of a function clone determine its topological structure? We pay particular attention to function clones which contain an oligomorphic permutation group, and discuss applications of this situation in model theory and theoretical computer science. 
\end{abstract}

\maketitle
\section{Introduction}
A \emph{function clone} (in the literature hitherto just \emph{clone}) 
over a set $D$
is a set of functions of finite arity on $D$
which is closed under composition and which contains the projections. Function clones appear naturally in algebra in the form of sets of term operations of algebras, which always form a function clone; indeed, every function clone is of this form. Since many important properties of an algebra, for example its subalgebras and its congruence relations, only depend on its term operations, function clones are of primordial importance in the understanding of algebras~\cite{KearnesKiss,HobbyMcKenzie}.
Function clones moreover generalize \emph{transformation monoids}, i.e., sets of \emph{unary} functions on a set $D$ closed under composition and containing the identity function. The latter generalize in turn \emph{permutation groups} on $D$, i.e., sets of permutations on $D$ closed under inverses and composition.

Similarly to (abstract) groups in group theory, \emph{abstract clones} have been studied extensively in universal algebra, though in disguise of \emph{varieties}~\cite{AbstractClones,KearnesKiss,HobbyMcKenzie}: roughly speaking, an abstract clone is an algebraic structure whose elements can be imagined as finitary functions on a fixed domain, together with composition operations on these elements and constant operations denoting the projections. Just as in the case of groups, every function clone gives rise to an abstract clone and vice-versa. Many insights about an algebra can be gained from the abstract clone associated with its term clone; this abstract clone basically encodes the equations which hold in the algebra~\cite{KearnesKiss,HobbyMcKenzie}.

Permutation groups carry a natural topology, the topology of pointwise convergence.
Under this topology, the corresponding group 
becomes a topological group since composition and taking inverses are continuous operations. 
Similarly, function clones are naturally equipped with the topology of pointwise convergence, and  again composition is continuous with respect to this topology.
As in the case of groups, where the study of \emph{topological groups} has without doubt been a fruitful venture for numerous fields of mathematics, 
it therefore makes sense to consider \emph{topological clones}, which consist of an abstract clone together with a topology on this structure under which the composition operations are continuous. The relationship between function clones and topological clones is in perfect analogy to the relationship between permutation groups and topological groups.  

Given the enormous literature on topological groups as the simultaneous topological and algebraic abstraction of permutation groups, and given the fact that the study of function clones and abstract clones constitute a considerable part of universal algebra, 
it is surprising that the analogous notion of topological clones has been entirely neglected in the literature. Inspired by the recent result~\cite{Topo-Birk} which indicates that a topological perspective on function clones in addition to the algebraic one is not only useful, but even inevitable if we strive for infinite versions of theorems about finite algebras, we here study for the first time topological clones explicitly. 
One of the purposes of this article is to demonstrate that topological clones exhibit a rich mathematical structure, and lead to many interesting and challenging problems. 
In particular, in this article we investigate the following research questions: 
\begin{itemize}
\item (Reconstruction of Topology) In which situations does the abstract clone of a function clone already determine its topology? 
\item  (Automatic Continuity) In which situations are homomorphisms or isomorphisms between function clones automatically continuous?
\end{itemize}

The corresponding questions for topological groups  have been the source of a wealth of beautiful results (see the survey article~\cite{Rosendal}),
in particular for automorphism groups of \emph{$\omega$-categorical} structures
(see, e.g.,~\cite{Barbina,Rubin,MacphersonBarbina,Herwig98,Truss,DixonNeumannThomas}). A countable structure $\Gamma$ is said to be $\omega$-categorical if and only if all countable models of the first-order theory of $\Gamma$ are isomorphic to $\Gamma$. We will be particularly interested in topological clones which arise in a similar
way from $\omega$-categorical structures, namely as their \emph{polymorphism clones}: for a structure $\Gamma$, the polymorphism clone $\pol(\Gamma)$ is the function clone consisting of all homomorphisms from finite powers of $\Gamma$ into  $\Gamma$.
Such topological clones have remarkable applications for which the answers to the questions posed above have strong consequences. We shall now discuss these applications.

\section{Three Applications}

\subsection{Reconstruction of $\omega$-categorical structures from their polymorphism clones}\label{sect:reconstructionofomega}
A permutation group is closed in the set of all permutations of its domain with respect to the topology of pointwise convergence if and only if it is the automorphism group $\Aut(\Gamma)$ of a relational structure $\Gamma$. It is natural to ask how much about a relational structure $\Gamma$ is coded in $\Aut(\Gamma)$. 
 If $\Gamma$ is $\omega$-categorical,
then $\Aut(\Gamma)$ determines $\Gamma$ up to \emph{first-order interdefinability}, that is, any structure $\Gamma'$ with 
$\Aut(\Gamma') = \Aut(\Gamma)$ has the property that all relations of $\Gamma$ have a first-order definition in $\Gamma'$ and vice-versa (see~\cite{Hodges}). 
In fact, this reconstruction property of a countably infinite structure is equivalent to $\omega$-categoricity. 
 
A function clone is closed in the set of all finitary functions on its domain with respect to the topology of pointwise convergence if and only if it is 
the polymorphism clone $\Pol(\Gamma)$ of a relational structure $\Gamma$.
The function clone $\Pol(\Gamma)$ 
encodes even more about $\Gamma$ than its automorphism group $\Aut(\Gamma)$. In particular, if $\Gamma$ is $\omega$-categorical and  $\Gamma'$ is such that $\Pol(\Gamma') = \Pol(\Gamma)$, then 
$\Gamma'$ and $\Gamma$ are \emph{primitive positive interdefinable},
that is, every relation in $\Gamma$ has a primitive positive definition  in $\Gamma'$ 
and vice versa~\cite{BodirskyNesetrilJLC}. 

On the other hand, if $\Gamma$ and $\Gamma'$ are countable $\omega$-categorical structures that only share the same automorphism group when it is viewed as a topological group rather than as a concrete permutation group, we obtain a different form of reconstruction~\cite{AhlbrandtZiegler}: the automorphism groups of $\Gamma$ and $\Gamma'$ are isomorphic as topological groups if and only if 
$\Gamma$ and $\Gamma'$ are \emph{first-order bi-interpretable} (see e.g.~\cite{Hodges}). 

It has been shown recently that the latter theorem and the theorem about polymorphism clones mentioned above can be naturally combined~\cite{Topo-Birk}: two countable $\omega$-categorical structures $\Gamma$ and $\Gamma'$ have isomorphic topological polymorphism clones if and only if $\Gamma$ and $\Gamma'$ are \emph{primitive positive bi-interpretable}. 
Figure~\ref{fig:reconstruction} gives a summary of all mentioned
forms of reconstruction of $\omega$-categorical structures.  

\begin{figure}
\begin{center}
\begin{tabular}{|lll|}
\hline
Reconstruction from: & Reconstruction up to: & Reference: \\ 
\hline
Permutation group &  First-order interdefinability & Ryll-Nardzewski~\cite{Hodges} \\
Topological group & First-order bi-interpretability & Ahlbrand and Ziegler~\cite{AhlbrandtZiegler} \\
Function clone & Primitive positive interdefinability & Bodirsky and Ne\v{s}et\v{r}il~\cite{BodirskyNesetrilJLC} \\
Topological clone & Primitive positive bi-interpretability & Bodirsky and Pinsker~\cite{Topo-Birk} \\
\hline
\end{tabular}
\end{center}
\caption{A schema for reconstruction of $\omega$-categorical structures.}
\label{fig:reconstruction}
\end{figure}

Positive answers to our two research questions combine nicely with the above result about primitive positive bi-interpretability:
when an $\omega$-categorical structure $\Gamma$ is such that 
clone isomorphisms between $\Pol(\Gamma)$ and other closed function clones are automatically homeomorphisms,
then this shows that already the \emph{abstract polymorphism clone} of $\Gamma$ determines $\Gamma$ up to primitive positive bi-interpretability. 

The analogous combination for groups has been studied intensively:
for many of the classical $\omega$-categorical structures it is known that they are determined by their abstract automorphism group up to first-order bi-interpretability. And indeed it is known to be consistent with ZF+DC that \emph{all} $\omega$-categorical structures are determined by their abstract automorphism group up to first-order bi-interpretability (\cite{Lascar}; cf.~the discussion in Section~8 in~\cite{Topo-Birk}).

\subsection{Complexity of Constraint Satisfaction Problems}
\label{ssect:csps}
Polymorphism clones and the topological clones they induce have applications in theoretical computer science. Every relational structure $\Gamma$ in a finite language defines a computational problem, called the \emph{constraint satisfaction problem} of $\Gamma$ and denoted by $\Csp(\Gamma)$, as follows: an instance of the problem is a primitive positive sentence $\phi$ in the language for $\Gamma$, i.e., a sentence of the form $\exists x_1,\dots,x_n (\phi_1 \wedge \dots \wedge \phi_m)$ where $\phi_1,\dots,\phi_m$ are atomic formulas; the problem is to decide whether or not $\phi$ holds in $\Gamma$. An instance of this problem therefore asks about the existence of elements of $\Gamma$ satisfying a given conjunction of atomic conditions. The structure $\Gamma$ is called the \emph{template} of the problem, and can be finite or infinite. Constraint satisfaction problems with infinite templates can model natural finite computational problems -- we refer to~\cite{Bodirsky-HDR, BodPin-Schaefer, BP-reductsRamsey, tcsps-journal} for an abundance of examples.

For finite and $\omega$-categorical structures $\Gamma$, the complexity of $\Csp(\Gamma)$ depends, up to polynomial-time interreducibility, only on $\Pol(\Gamma)$. More precisely,  if $\Gamma$ and $\Gamma'$ are $\omega$-categorical structures in finite relational languages on the same domain, and if $\pol(\Gamma')=\pol(\Gamma)$, then $\Csp(\Gamma)$ and $\Csp(\Gamma')$ are polynomial-time equivalent 
 (cf.~\cite{JBK,JeavonsAlgebra, BodirskyNesetrilJLC}). This fact is the basis of what is known as the \emph{algebraic approach} to constraint satisfaction. But the algebraic approach goes even further: for finite structures $\Gamma$ the complexity of $\Csp(\Gamma)$  depends up to polynomial time only on $\Pol(\Gamma)$, viewed as an abstract clone~\cite{JBK,JeavonsAlgebra}. In the $\omega$-categorical case, it has been shown recently to depend only on $\Pol(\Gamma)$, viewed as a topological clone~\cite{Topo-Birk}. Moreover, up to now no two $\omega$-categorical structures with abstractly isomorphic polymorphism clones but CSPs of different (up to polynomial-time reductions) complexity are known, and it has been shown recently
that at least some aspects of the CSP of an $\omega$-categorical structure are captured by the algebraic structure of
$\Pol(\Gamma)$~\cite{BartoPinsker}.
 
\subsection{Pseudovarieties of oligomorphic algebras}\label{sect:pseudovarieties}

The \emph{term clone} $\Clo(\mathfrak A)$ of an algebra $\frak A$ with signature $\tau$ is the set of all functions with finite arity on the domain of $\frak A$ which can be written as $\tau$-terms over $\frak A$. Clearly, $\Clo(\mathfrak A)$ is always a function clone, and all function clones are of this form.

Let ${\frak A}$, ${\frak B}$ be algebras of the same signature $\tau$. The assignment which sends every term function over $\frak A$ to the corresponding term function over $\frak B$ is a well-defined function from $\Clo(\mathfrak A)$ to $\Clo(\mathfrak B)$ if and only if all equations which hold between terms over $\frak A$ also hold over $\frak B$. In that case, it is in fact a surjective \emph{clone homomorphism}, i.e., it preserves projections and composition of functions (cf.~Section~\ref{sect:prelims}); it is then called the  \emph{natural homomorphism} from $\Clo(\mathfrak A)$ onto $\Clo(\mathfrak B)$.

A \emph{pseudovariety} is a class of algebras of the same signature which is closed under subalgebras, homomorphic images, and finite products. Since these operators are among the most fundamental and natural for algebras, pseudovarieties play an important role in the study of algebras. The pseudovariety \emph{generated} by a \emph{finite} algebra ${\frak A}$, i.e., the smallest pseudovariety which contains ${\frak A}$, was characterized in a classical theorem due to Garrett Birkhoff via $\Clo(\mathfrak A)$, viewed as an abstract clone: it contains precisely those finite algebras $\frak B$ for which the natural homomorphism from $\Clo(\mathfrak A)$ onto $\Clo(\mathfrak B)$ exists (\cite{Bir-On-the-structure}; cf.~ also  Exercise 11.5 in combination with the proof of Lemma 11.8 in~\cite{BS}).

Birkhoff's theorem has recently been generalized to \emph{oligomorphic algebras}. A permutation group on a countable set $D$ is called \emph{oligomorphic} iff its componentwise action on any finite power of $D$ has finitely many orbits. A function clone is called oligomorphic iff it contains an oligomorphic permutation group. It follows from the theorem of Ryll-Nardzewski (see~\cite{Hodges}) that the closed oligomorphic clones are precisely the polymorphism clones of $\omega$-categorical structures. An algebra is oligomorphic iff the topological closure of its term clone is oligomorphic and hence the polymorphism clone of an $\omega$-categorical structure. 

It is easy to see that all elements of the pseudovariety generated by an oligomorphic algebra $\frak A$ must be finite or oligomorphic. Now the generalization of Birkhoff's theorem states that if ${\frak B}$ is an oligomorphic or finite algebra in the signature of $\frak A$, then ${\frak B}$ is contained in the pseudovariety generated by ${\frak A}$ if and only if the natural homomorphism from the closure of $\Clo(\frak A)$ to the closure of $\Clo(\frak B)$ exists and is continuous~\cite{Topo-Birk}. In Birkhoff's finite version of this theorem, there is of course no continuity condition on the natural homomorphism, since function clones on a finite domain are discrete and so homomorphisms from finite function clones are always continuous. The present paper addresses the question for which oligomorphic algebras we can drop continuity in the generalized theorem.

\section{Main Notions, More Background, and Results}
\label{sect:prelims}
We introduce the notion of a topological clone, and recall the definitions of a function clone and abstract clone in more detail. We then define variants of several reconstruction notions from the literature on topological groups for topological clones, and give an overview of the results we will obtain.

\subsection{Function clones and abstract clones}

\begin{defn}
A \emph{function clone} $\mathscr C$ (in the literature simply \emph{clone}) over a set $D$ is a set of functions of finite arity over $D$ such that 
\begin{itemize}
\item $\mathscr C$ contains for all $1\leq k\leq n<\omega$ the \emph{$k$-th $n$-ary projection $\pi^n_k\colon D^n\To D$}, uniquely defined by the equation $\pi^n_k(x_1,\ldots,x_n)=x_k$;
\item whenever $f\in \fC$ is $n$-ary, and $g_1,\ldots,g_n\in\fC$ are $m$-ary, then the $m$-ary function $f(g_1,\ldots,g_n)$ defined by
$$(x_1,\ldots,x_m)\mapsto f(g_1(x_1,\ldots,x_m),\ldots,g_n(x_1,\ldots,x_m))$$ is  an element of $\fC$.
\end{itemize}
We write $\mathscr C^{(n)}$ for the $n$-ary functions in $\mathscr C$, for all $n\geq 1$. The set $D$ is also called the \emph{domain} of $\mathscr C$, and the elements of $\mathscr C$ are also referred to as the \emph{operations} of $\mathscr C$.
\end{defn}

An important source of examples for function clones
are \emph{polymorphism clones} of structures.
For a structure $\Gamma$ with domain $D$, a \emph{polymorphism} of $\Gamma$ is a homomorphism from $\Gamma^n$ to $\Gamma$ for some $n\geq 1$. It is easy to verify that the set 
$\Pol(\Gamma)$ of all polymorphisms of a structure $\Gamma$ is a function clone.

The algebraic structure of function clones can best be understood via the appropriate notion of a clone homomorphism.

\begin{defn} Let $\fC, \fD$ be function clones (not necessarily over the same set). Then a function $\xi \colon \fC\To\fD$ is called a \emph{(clone) homomorphism} iff 
\begin{itemize}
\item it preserves arities of functions;
\item for all $1\leq k\leq n<\omega$, the $k$-th $n$-ary projection in $\fC$ is sent to the $k$-th $n$-ary projection in $\fD$;
\item $\xi(f(g_1,\ldots,g_n))=\xi(f)(\xi(g_1),\ldots,\xi(g_n))$ whenever $n,m\geq 1$, $f\in \fC$ is $n$-ary, and $g_1,\ldots,g_n\in\fC$ are $m$-ary.
\end{itemize}
\end{defn}
 
The structure given by this notion of a clone homomorphism can also be formalized abstractly; and where from permutation groups we obtain (abstract) groups, we obtain \emph{(abstract) clones} from function clones. In practice we will not need this formalization, but only the corresponding notion of homomorphism, but we include the definition for the sake of completeness.

\begin{defn}
A \emph{clone} ${\mathfrak C}$ (in the literature \emph{abstract clone}) is a multi-sorted structure with sorts $\{C^{(n)} \; | \; n \geq 1\}$ and the
signature $\{\pi_k^n \; | \; 1\leq k \leq n\}
\cup \{\comp^n_m \; | \; n,m \geq 1\}$. The elements of the sort $C^{(n)}$, which we also denote by ${\mathfrak C}^{(n)}$, will be called the \emph{$n$-ary operations} of ${\mathfrak C}$. 
We denote a clone by 
$${\mathfrak C} = (C^{(1)},C^{(2)},\dots;(\pi_k^n)_{1\leq k \leq n},(\comp^n_m)_{n,m \geq 1})$$ and require that
$\pi_k^n$ is a constant in $C^{(n)}$,
and that  
$\comp_m^n \colon C^{(n)} \times (C^{(m)})^n \to C^{(m)}$ is an
 operation of arity $n+1$.
Moreover, 
$\comp_n^n(f,\pi_1^n,\dots,\pi_n^n) = f$, $\comp^n_m(\pi_k^n,f_1,\dots,f_n) = f_k$, and 
\begin{align*}
& \comp_m^n(f,\comp_k^m(g_1,h_1,\dots,h_m),\dots,\comp_k^m(g_n,h_1,\dots,h_m))  \\
= & \comp^m_k(\comp^n_m(f,g_1,\dots,g_n),h_1,\dots,h_m) \; .
\end{align*}
\end{defn}

We also write $f \circ (g_1,\dots,g_n)$ or $f(g_1,\dots,g_n)$ instead
of $\comp^n_m(f,g_1,\dots,g_n)$ when 
$m$ is clear from the context. When composing unary elements $f$ and $g$, we also write $fg$ for better readability.

There is a more
convenient way to write equations that hold in a
clone; e.g., we will say that $f \in C^{(2)}$ \emph{satisfies} $\forall x,y. f(x,y)=f(y,x)$, or $f(x,y)=f(y,x)$ \emph{holds in $\mathfrak C$} instead of writing $\comp^2_2(f,\pi_1^2,\pi_2^2)=\comp^2_2(f,\pi_2^2,\pi_1^2)$; this can be viewed as syntactic sugar. Note that a homomorphism from a clone $\mathfrak C$ to a clone $ \mathfrak D$ is just a function which preserves all equations that hold in $\mathfrak C$. An example of an equation which will be important throughout the paper is the following. 
A unary element $e$ of a clone is called \emph{invertible} iff there exists a unary element $f$ in the clone
such that $\forall x. f(e(x))=e(f(x))=x$ is satisfied. Clearly the invertible elements of a clone form an abstract group, and the unary elements of a clone form an abstract monoid with the composition operation $\comp^1_1$.

Every function clone $\mathscr C$ gives rise to an (abstract) 
clone $\mathfrak C$ in the obvious way.
Conversely, a straightforward generalization of Cayley's theorem
for groups shows that for every clone $\mathfrak C$ 
there exists a function clone whose abstract clone is $\mathfrak C$. We call any such realization of $\mathfrak C$ as a function clone $\mathscr C$ on a set $D$ an \emph{action} of $\mathfrak C$ on the set $D$.

\subsection{Topological clones}
\label{sect:topo-clones}
On any set $D$, there is a largest function clone $\mathscr O_D$, which consists of all finitary operations on $D$. The set $\mathscr O_D$ is naturally equipped with the topology of pointwise convergence, with respect to which the composition of functions is continuous. A basis of open sets of this topology 
is given by the sets of the form
$$\{ f \colon D^n\To D \; | \; 
f(a^1_1,\dots,a^1_n) = a^1_0, \dots,f(a^m_1,\dots,a^m_n) = a_0^m \} \;.$$
For countably infinite $D$, $\mathscr O_D$ becomes a Polish space with this topology; in fact, $\mathscr O_D$ is then homeomorphic to the Baire space. A compatible complete metric can be defined as follows. 
For each $n$, we 
fix an enumeration $a_1^n,a_2^n,\dots$ of $D^n$.
When $f,g\in\mathscr O_D$ have the same arity $n$, 
then put $d(f,g) = 1/2^{\min(i \; | \; f(a_i^n) \neq g(a_i^n))}$.
When $f$ and $g$ have distinct arity, 
put $d(f,g) = 1$.

The function clones on $D$ which are closed in $\mathscr O_D$ with respect to this topology are precisely the clones of the form $\Pol(\Gamma)$ for some first-order structure $\Gamma$ with domain $D$. As a subset of $\mathscr O_D$, any function clone on $D$ inherits a topology from $\mathscr O_D$.  Hence, it carries a topological structure in addition to its algebraic structure, motivating the following new definition.

\begin{defn}
A \emph{topological clone} $\mathbf C$ is a clone $${\mathfrak C} = (C^{(1)},C^{(2)},\dots;(\pi_k^n)_{1\leq k \leq n},(\comp^n_m)_{n,m \geq 1})$$ together with a topology on $\bigcup_{n\geq 1}C^{(n)}$ such that each $C^{(n)}$ is a clopen set and such that the composition operations are continuous.
\end{defn}

As discussed above, every function clone gives rise to a topological clone. We will be interested in topological clones induced by function clones on a countably infinite set $D$, and
write $\mathbf O$ 
for the topological clone of the function clone $\mathscr O_D$; cf.~\cite{GoldsternPinsker} for a survey of function clones on $D$. Moreover, we write $\mathbf S$ for the topological group of the full symmetric group $\mathscr S_D$ over a countably infinite set $D$. It is known that the closed subgroups of ${\mathbf S}$ 
are precisely those topological groups that are Polish
and have a left-invariant ultrametric~\cite{BeckerKechris}. There exists an analogous characterization of the closed subclones of $\mathbf O$, obtained after the completion of the first version of the present article~\cite{BodirskySchneiderTopological}; confer the open problems section (Section~\ref{sect:open}).

\subsection{Topological monoids}

We are going to consider various reconstruction notions for topological groups and clones; in particular, we  will use known reconstruction results for groups to obtain such results for clones. A natural class of objects between the two classes is the class of monoids. Here we distinguish \emph{transformation monoids}, i.e., sets of unary functions on a fixed set which are closed under composition and which contain the identity function; \emph{(abstract) monoids}, with their well-known definition; and \emph{topological monoids}, i.e., abstract monoids which in addition carry a topology under which composition is continuous. We denote the transformation monoid of all unary functions on a set $D$ by $\mathscr O^{(1)}_D$, and write $\mathbf O^{(1)}$ for topological monoid induced by $\mathscr O^{(1)}_D$ when $D$ is countably infinite. For any set $D$, the closed submonoids of $\mathscr O^{(1)}_D$ are precisely the endomorphism monoids of first-order structures on $D$; if $\Gamma$ is such a structure, then we write $\End(\Gamma)$ for its endomorphism monoid.

Every permutation group (abstract group, topological group) can be seen as a transformation monoid (abstract monoid, topological monoid); conversely, the invertible elements of a transformation monoid (abstract monoid, topological monoid) form a permutation group (abstract group, topological group).  We would like to point out that $\mathbf S$ is not closed in $\mathbf O^{(1)}$: whenever $f\in \mathscr O^{(1)}_D$ is injective, where $D$ is countable, then on every finite subset of $D$ there exists a permutation on $D$ which agrees with $f$ on this finite subset, and hence $f$ is an element of the closure of $\mathscr S_D$ . Consequently, closed subgroups of $\mathbf S$ need not be closed in $\mathbf O^{(1)}$.

Similarly, monoids can be interpreted as clones by adding the projections and closing under composition; and conversely, the set of unary functions of a function clone (or the unary elements of an abstract clone) form a transformation monoid (an abstract monoid). Some of our examples of topological clones will really be examples of topological monoids. 

It is worth noting, however, that not every monoid homomorphism can be extended to a clone homomorphism between the corresponding clones; indeed, there is a slight technical condition which has to be added in order to ensure this. An $n$-ary element $f$ of a clone is called a \emph{constant} iff $\forall x_1,\ldots,x_n,y_1,\ldots,y_n.\; f(x_1,\ldots,x_n)=f(y_1,\ldots,y_n)$ holds in the clone; note that this formula can be written without quantifiers in the language of clones. In the language of monoids, constants cannot be defined without quantifiers, a fact which is reflected in the following proposition.

\begin{proposition}\label{prop:monoids-vs-clones}
Let $\frak M$, $\frak N$ be monoids, and let $\frak M'$, $\frak N'$ be the corresponding clones. Then:
\begin{itemize}
\item the restriction of any homomorphism $\xi\colon \frak M'\To \frak N'$ to $\frak M$ is a monoid homomorphism;
\item the natural extension of a homomorphism $\xi\colon \frak M\To \frak N$ to $\frak M'$ is a clone homomorphism if and only if $\xi$ sends constants (in $\frak M'$) to constants (in $\frak N'$).
\end{itemize}
\end{proposition}
\begin{proof}
The first statement is clear, since every equation in the language of monoids is also an equation in the language of clones. Now consider the second item. If $\xi$ sends a constant to a non-constant, then its natural extension is not a clone homomorphism. For the converse, it is easily verified that all equations which hold in $\frak M'$ are equivalent to equations of the form $\forall x,y.\; f(x)=g(y)$ or $\forall x.\; f(x)=g(x)$. The first type is equivalent to the two equations $\forall x,y.\; f(x)=f(y)$ and $\forall x.\; f(x)=g(x)$. These two are preserved by $\xi$: the first one because $\xi$ preserves constants, and the second one is trivially preserved because $\xi$ is a function.
\end{proof}

Observe that as for continuity questions, it is irrelevant whether we see a topological monoid as a monoid or as a clone; that is, a homomorphism between topological monoids is continuous iff its (unique) extension to the corresponding topological clones is.

The link between topological clones and primitive positive bi-interpretability mentioned in Section~\ref{sect:reconstructionofomega} has an analog for topological monoids and existential positive bi-interpretability, which however
does not hold for all $\omega$-categorical structures~\cite{BodJunker}: the reason for this is the inability of monoids to express constants exhibited in Proposition~\ref{prop:monoids-vs-clones}. 

Our motivation to also study reconstruction for monoids is two-fold: firstly, several challenges for clones already become apparent in the
less complex case of monoids; secondly, our reconstruction results for clones typically build on 
reconstruction results for the monoids given by their unary functions. 

\subsection{Oligomorphicity}

We will be interested in topological clones induced by oligomorphic function clones on a countable set (recall the definition of \emph{oligomorphic} in Section~\ref{sect:pseudovarieties}). Whether or not a topological group $\mathbf G$ has a continuous action that induces an oligomorphic permutation group has an elegant
characterization using the terminology from Polish groups (without referring to any particular action of $\mathbf G$); confer~\cite{Tsankov}. 
It also turns out that if $\mathbf G$
has such an action, then \emph{every} continuous action of $\mathbf G$ on a countable set with finitely
many orbits induces an oligomorphic permutation group.
We therefore call such topological groups \emph{oligomorphic}. 
We say a subclone of $\mathbf O$ (submonoid of $\mathbf O^{(1)}$) is \emph{oligomorphic}
iff the set of invertible elements of the clone (monoid) forms an oligomorphic
topological group.

\subsection{Reconstruction notions}

We study the question whether we can reconstruct the topology of closed subclones of the topological clone $\mathbf O$ (or equivalently, of the function clone $\mathscr O_D$) from the abstract clone structure alone. In this context, several reconstruction notions make sense. The following definitions are inspired by the literature on topological groups.

\begin{defn}
Let $\mathbf C$ be a closed subclone of $\mathbf O$. We say that
\begin{itemize}
\item $\mathbf C$ is \emph{reconstructible} (or that $\mathbf C$ has \emph{reconstruction}) iff for every other closed subclone $\mathbf D$ 
of $\mathbf O$, if there exists a clone isomorphism between $\mathbf C$ and $\mathbf D$, then there also exists a clone isomorphism between $\mathbf C$ and $\mathbf D$ which is a homeomorphism;
\item $\mathbf C$ has \emph{automatic homeomorphicity} iff every clone isomorphism between $\mathbf C$ and a closed subclone of $\mathbf O$ is a homeomorphism;
\item $\mathbf C$ has  \emph{automatic continuity} iff every clone homomorphism from $\mathbf C$ into $\mathbf O$ is continuous.
\end{itemize}
All these notions are analogously defined for closed subgroups of $\mathbf S$ and closed submonoids of $\mathbf O^{(1)}$. 
\end{defn}

Note that automatic homeomorphicity implies reconstruction; otherwise, the notions are not obviously related. However, for groups automatic continuity implies automatic homeomorphicity.
 
\begin{proposition}[Corollary 2.8 in~\cite{Lascar}]
\label{prop:homeo-from-cont}
Any continuous isomorphism between closed subgroups of ${\mathbf S}$ is a homeomorphism. 
\end{proposition}

Proposition~\ref{prop:homeo-from-cont} shows that automatic continuity of closed subgroups of ${\mathbf S}$ is 
a property of the abstract group in the sense that
if two closed subgroups of ${\mathbf S}$ are isomorphic as abstract groups, and one has automatic continuity,
then so has the other. 

\subsection{The situation for topological groups}
\label{sect:groups}

There are two dominant methods for proving reconstruction of group topology; the two methods have in common that they imply  reconstruction via automatic homeomorphicity. A method for proving reconstruction which would apply also to groups without automatic homeomorphicity seems hardly conceivable.

The first method is showing \emph{automatic continuity}, or equivalently, the \emph{small index property}.  A topological group $\G$ has the small index property iff every subgroup of $\G$ of at most countable index is open. It is a folklore fact, and not difficult to show, that a topological group has automatic continuity if and only if it has the small index property. The small index property has been verified for the following groups:

\begin{itemize}
\item $\mathbf S$~\cite{Rabinovich,Semmes,DixonNeumannThomas} -- that is, $\Aut({\mathbb N};=)$;  
\item the automorphism groups of countable vector spaces over finite fields~\cite{EvansSmallIndex};
\item $\Aut({\mathbb Q};<)$ and the automorphism group of the atomless Boolean algebra~\cite{Truss}; 
\item the $\omega$-categorical
dense semi-linear order giving rise to a meet-semilattice~\cite{DrosteHollandMacphersonII};
\item the automorphism group of the random graph~\cite{HodgesHodkinsonLascarShelah};
\item all automorphism groups of $\omega$-categorical $\omega$-stable structures~\cite{HodgesHodkinsonLascarShelah};
\item the automorphism groups of the Henson graphs~\cite{Herwig98}.
\end{itemize}

The second method for proving reconstruction
is Rubin's \emph{forall-exists interpretations}. 
Their attractive feature is that they allow
us to recover an automorphism group as a permutation group from the abstract group structure when one restricts the category of structures whose automorphism groups one is  interested in.
 More precisely, if $\Gamma$ is an $\omega$-categorical structure which has weak forall-exists interpretations,
and if $\Delta$ is another $\omega$-categorical structure without algebraicity and with an isomorphism $\xi$ between $\Aut(\Gamma)$ and $\Aut(\Delta)$, then there exists a bijection $i$ between the domain of $\Gamma$ and the domain of $\Delta$ such that 
for all $\alpha \in \Aut(\Gamma)$ we have
$\forall x. \, \xi(\alpha)(x) = i(\alpha(i^{-1}(x)))$. Discussing forall-exists interpretations is beyond the scope of this paper, but they have been given for:

\begin{itemize}
\item the random graph, $({\mathbb Q};<)$, all homogeneous countable graphs, and various $\omega$-categorical semi-linear orders~\cite{Rubin};
\item the universal homogeneous partial ordering, the universal homogeneous tournament~\cite{Rubin} (for those structures it is not known whether they have the small index property);
\item universal homogeneous $k$-hypergraphs, and the Henson digraphs~\cite{MacphersonBarbina}. 
\end{itemize}

\subsection{Density of the invertibles} 
It is a fact that every CSP of an $\omega$-categorical structure $\Gamma$ is equal (as the set of instances with a positive answer) 
 to the CSP of another $\omega$-categorical structure $\Gamma'$ which is a \emph{model-complete core}~\cite{BodHilsMartin}, i.e., a structure whose automorphisms are dense in its endomorphisms; in other words, $\Aut(\Gamma')$ is dense in the unary part of $\pol(\Gamma')$. We are therefore particularly interested in this situation. Surprisingly, it turns out to be a non-trivial task to show for a given closed oligomorphic subgroup $\mathbf G$ of $\mathbf S$ with automatic continuity (the strongest form of reconstruction) that the closure $\overline{\mathbf G}$ of $\mathbf G$ in $\mathbf O^{(1)}$ has some form of reconstruction. Note that when $\Gamma$ is a structure whose automorphism group induces $\mathbf G$, then $\overline{\mathbf G}$ is the topological monoid induced by the monoid $\overline{\Aut(\Gamma)}$ of elementary self-embeddings of $\Gamma$, i.e., of self-embeddings preserving all first-order formulas over $\Gamma$.
 
 The related problem of proving that oligomorphic groups are isomorphic as topological groups when their closures are isomorphic as abstract (non-topological) monoids has been investigated in~\cite{Lascar}; confer also Section~\ref{sect:monoids} for a result we will use from this work.
 


\subsection{Results}

\subsubsection{Positive results} 

We present methods for proving automatic homeomorphicity of closed subclones $\mathbf C$
of $\mathbf O$. Our first result concerns the closure in $\mathbf O^{(1)}$ of $\mathbf S$.

\begin{itemize}
\item If $\Delta$ is a homogeneous structure over a finite relational language without algebraicity, with the \emph{joint extension property}, and such that $\Aut(\Delta)$ has automatic continuity, then its closure $\overline{\Aut(\Delta)}$ has automatic homeomorphicity (Section~\ref{sect:monoids}).
\end{itemize}

The task of proving automatic homeomorphicity for other clones is then split into
proving that isomorphisms between $\bf C$ and
closed subclones of $\bf O$ are continuous,
and proving that these isomorphisms are open. For proving continuity, we present a technique based on so-called \emph{gates} as well as the above result: 
\begin{itemize}
\item If $\mathbf M$ is a closed submonoid of $\mathbf O^{(1)}$ such that the closure of the set of its invertible elements has automatic homeomorphicity, and which has a \emph{gate} with respect to this closure, then every isomorphism from $\mathbf M$ onto another closed submonoid of $\mathbf O^{(1)}$ is continuous (Section~\ref{sect:gatesforunary}).
\item If $\mathbf C$ is a closed subclone of $\mathbf O$ which has a \emph{gate covering
}, and if $\xi$ is an isomorphism from $\mathbf C$ onto another closed subclone of $\mathbf O^{(1)}$ whose restriction to the unary elements is continuous, then $\xi$ is continuous (Section~\ref{subsect:gatescontinuity}).
\end{itemize}
We also present another technically unrelated method for proving continuity via Birkhoff's HSP theorem~\cite{Bir-On-the-structure} in  Section~\ref{subsect:birkhoffcont}.

Concerning openness, we first obtain two results serving different classes of clones: 
\begin{itemize}
\item function clones that contain all constant functions in Section~\ref{subsect:constantsopenness};
\item \emph{transitive} function clones in Section~\ref{subsect:transitivityopenness}.
\end{itemize}
We then show how a recent topological variant of Birkhoff's theorem from~\cite{Topo-Birk} can be exploited in this context in Section~\ref{subsect:topobirkopenness}.

Using these general results and methods, we obtain automatic homeomorphicity for several transformation monoids and function clones.
\begin{itemize}
\item The monoids of self-embeddings of the empty structure, the random graph, and the random tournament have automatic homeomorphicity (Section~\ref{sect:monoids}).
\item The \emph{Horn clone} and the polymorphism clone of the random graph 
have automatic homeomorphicity (Section~\ref{sect:clones}).
\item Any closed subclone of $\mathscr O_\omega$ containing $\mathscr O^{(1)}_\omega$ has automatic continuity and automatic homeomorphicity (Section~\ref{subsect:birkhoffcont}).
\end{itemize}

\subsubsection{Negative results}

On the negative side, we will show the following.

\begin{itemize}
\item There exists a closed oligomorphic submonoid $\mathbf M$ of $\mathbf O^{(1)}$ and an isomorphism $\xi\colon \mathbf M\To \mathbf M$ which is not continuous; in particular, $\mathbf M$ does not have automatic homeomorphi-city (Section~\ref{sect:monoids}). The example lifts to clones via Proposition~\ref{prop:monoids-vs-clones}.
\item There are simple conditions on monoids which imply that they cannot have automatic continuity; in particular, no monoid of self-embeddings of an $\omega$-categorical structure has automatic continuity (Section~\ref{sect:monoids}).
\end{itemize}
The second statement stands in sharp contrast with the situation for groups, and indicates that the notion of automatic continuity is somewhat too strong for topological monoids.

\section{Topological Monoids}
\label{sect:monoids}

\subsection{Negative results}
\begin{theorem}\label{thm:no-auto-homeo}
There exists a closed oligomorphic submonoid $\mathbf M$ of $\mathbf O^{(1)}$ which has an automorphism $\xi\colon \mathbf M\To \mathbf M$ that is not continuous. In particular, $\mathbf M$ does not have automatic homeomorphicity. Moreover, $\xi$ sends constants to constants, and hence $\xi$ lifts to a discontinuous automorphism of the corresponding clone.
\end{theorem}

\begin{proof}
Let $(S^n)_{n\geq 2}$ be a sequence of relational symbols such that $S^n$ is $n$-ary for all $n\geq 2$. Consider the class of all finite structures in this language such that $S^n$ is interpreted as a totally symmetric $n$-ary relation of injective tuples, and let $(V_0;(S_0^n)_{n\geq 2})$ and $(V_1;(S_1^n)_{n\geq 2})$ be two copies of the Fra\"{i}ss\'{e} limit of this class with disjoint domains $V_0$ and $V_1$. Fix an isomorphism $\iota \colon V_0 \rightarrow V_1$ between them.
 Set $V:=V_0\cup V_1$ and $\alpha:=\iota\cup\iota^{-1}$; then $\alpha$ is a permutation on $V$ which is equal to its own inverse. For notational simplicity, we write $E_i:=S^2_i$, for $i\in\{0,1\}$; clearly, $(V_i,E_i)$ is isomorphic to the random graph.  
Let $E:=E_0\cup E_1\cup \{(v,\alpha(v))\; |\; v\in V \}$, and $S^n_{01}:=S^n_0\cup S^n_1$ for all $n\geq 3$.
For any self-embedding $e$ of $(V_0;E_0)$, let $\bar{e}$ be the self-embedding of $(V;E)$ defined by $\bar e(v):=e(v)$ if $v\in V_0$, and $\bar e(v):=\alpha(e(\alpha(x)))$ if $v\in V_1$. Note that embeddings of the form $\bar{e}$ commute with $\alpha$, i.e., $\bar e\circ\alpha=\alpha\circ \bar e$.\smallskip

\noindent {\it Claim.} Let $f$ be a self-embedding of $(V;E)$. Then either there exists a self-embedding $e$ of $(V_0;E_0)$ such that $f=\bar e$ or $f=\alpha\bar{e}=\bar{e}\alpha$, or the range of $f$ is contained in $V_i$ for some $i\in \{0,1\}$.\smallskip

To verify the claim, take any $u\in V_0$, and assume that $f(u)\in V_i$, where $i\in\{0,1\}$. 
We show that the neighbors of $u$ in $(V_0;E_0)$ are also mapped to $V_i$ under $f$. So let $v\in V_0$ be such that $(u,v)\in E_0$. If $f(v)=\alpha(f(u))$, then let $w\in V_0$ be so that $\{u,v,w\}$ induces a complete graph in $(V_0;E_0)$. Then $f(w)$ must be connected to $f(u)$ and $f(v)$ in $(V;E)$, a contradiction. Hence, $f(v)\neq \alpha(f(u))$. Thus $f(v)\in V_i$, as $f(u)$ is in $(V;E)$ adjacent to no element of $V_{1-i}$  except for $\alpha(f(u))$. As the random graph has diameter 2, we obtain that $f[V_0]$ is contained in either $V_0$ or $V_1$.  Similarly, $f[V_1]$ is contained in either $V_0$ or $V_1$. Hence there are four possibilities.  If $f$ maps $V$ into $V_0$ or $V_1$, then we are done. 
 If $f[V_i]\subseteq V_i$ for both $i\in \{0,1\}$, then $f$ is of the form $\bar{e}$. If $f[V_0]\subseteq V_1$ and $f[V_1]\subseteq V_0$, then $\alpha f$ is as in the preceding case, and hence $\alpha f=\bar{e}$ for a self-embedding $e$ of $(V_0;E_0)$. Consequently,  $f=\alpha\bar{e}$.\smallskip

The claim implies that 
$\alpha$ commutes with all elements of $\overline{\Aut(V;E)}$. However, whenever the range of a self-embedding of $(V;E)$ is contained in $V_0$, then this embedding does not commute with $\alpha$.
\smallskip

Now take the structure $(V;E,(S^n_{01})_{n\geq 3})$, and consider the structure $\Delta$ which consists of two copies $\Gamma$, $\Gamma'$ of this structure on disjoint domains $V$ and $V'$, plus an extra element $c$ outside $V\cup V'$ and predicates for $V$ and $V'$. Write $D$ for the domain $V\cup V'\cup\{c\}$ of $\Delta$; abusing the notation slightly, we write $E$ and $(S^n_{01})_{n\geq 3}$ for the relations of $\Delta$ (so each of these relations is the union of the relations with the same name in $\Gamma$ and in $\Gamma'$). Thus the automorphism group of $\Delta$ really is $\Aut(\Gamma)\times \Aut{(\Gamma')}$. In particular, it is oligomorphic since $\Aut(\Gamma)$ is oligomorphic: this follows readily from the easily verified fact that $(V;E,(S^n_{01})_{n\geq 3},V_0,V_1)$ is homogeneous.

In the following, we say that a function $f\colon D\To D$ \emph{eradicates} a relation $S^n_{01}$ iff $S^n_{01}$ holds for no tuple in the range of $f$. Now writing $\Emb(\Xi)$ for the monoid of self-embeddings of a structure $\Xi$, set for all $3\leq m\leq \omega$
\begin{align*}
\mathscr F_m:=\{f\in \Emb(D;E,(S^n_{01})_{n\geq m})\;|\; 
&f[V]\subseteq V',\;f[V'\cup\{c\}]= \{c\}, \text{ and }\\ 
&f\text{ eradicates all } S^n_{01}\text{ with }3\leq n<m\}.
\end{align*}
Writing $f_c\colon D\To D$ for the constant function with value $c$, set moreover
\begin{align*}
\M_\infty&:=\overline{\Aut(\Delta)}\cup\{f_c\}\cup \mathscr F_\omega,\\
\M_{<\infty}&:=\overline{\Aut(\Delta)}\cup\{f_c\}\cup \bigcup_{3\leq m<\omega} \mathscr F_m,\text{ and }\\
\M&:=\M_\infty\cup \M_{<\infty }.
\end{align*}

Note that $\M$ is a closed monoid, that 
$\M_\infty$ is a closed submonoid of $\M$, and that $\M_{<\infty }$ is a submonoid of $\M$ which is dense in $\M$.  Moreover, $\M_\infty\cap \M_{<\infty }=\overline{\Aut(\Delta)}\cup\{f_c\}$.

Recall the function $\alpha$, and imagine it acts on $D$ keeping $c$ fixed, and acting on $V$ and $V'$ as above. Now define a mapping $\xi\colon \M\To \M$ by
$$
\xi(f):=\begin{cases}f &, f\in \M_{<\infty} \\\alpha\circ f\circ \alpha&, f\in\M_\infty.
\end{cases}
$$
Then $\xi$ is well-defined because $\alpha$ commutes with all functions in $\M_\infty\cap \M_{<\infty }=\overline{\Aut(\Delta)}\cup\{f_c\}$. 
Clearly, the restriction of $\xi$ to $\M_\infty$ and $\M_{<\infty }$, respectively, is an inner automorphism of those monoids. 

We claim that $\xi$ is an automorphism of $\M$. To see this, note first that we have already observed that its restriction to $\M_\infty$ and $\M_{<\infty }$, respectively, is an automorphism. So let $f\in \M\setminus \M_{<\infty}$ and $g\in \M\setminus \M_{\infty}$ be given. Then $f\circ g=f_c=\xi(f\circ g)$ and $\xi(f)\circ \xi(g)=f_c$, proving $\xi(f\circ g)= \xi(f)\circ \xi(g)$. Similarly, $\xi(g\circ f)= f_c = \xi(g)\circ \xi(f)$.

However, if $g\in\M_{\infty}$ is so that it does not commute with $\alpha$, and $(f_n)_{n\in\omega}$ is a sequence in $\M_{<\infty}$ converging to $g$, then $(\xi(f_n))_{n\in\omega}$ will still converge to $g$, proving that $\xi$ is not continuous.

\end{proof}

We will now see that many closed submonoids of $\mathbf O^{(1)}$ do not have automatic continuity, so that this notion is arguably less useful than for closed subgroups of $\mathbf S$.

\begin{proposition}\label{prop:no-automatic-cont}
Let ${\mathbf M}$ be a closed submonoid of 
$\mathbf O^{(1)}$. Suppose that ${\mathbf M}$ contains a submonoid  $\mathbf N$ such that 
\begin{enumerate}
\item $\mathbf N$  is not closed in ${\mathbf M}$;
\item composing any element of $\mathbf M$ with an element outside $\mathbf N$ yields an element outside $\mathbf N$.
\end{enumerate}
Then ${\mathbf M}$ does not have automatic continuity.
\end{proposition}
\begin{proof}
Let $D$ be a countable set, and let $\M\subseteq \mathscr O^{(1)}_D$ be an action of $\mathbf M$ on $D$. Write $\mathscr N$ for the transformation monoid corresponding to $\mathbf N$. Let $i$ be a bijection between $D$ and $D \setminus \{c\}$ for
some $c \in D$. We define a monoid homomorphism $\xi$
from ${\mathscr M}$ to ${\mathscr O^{(1)}_D}$. Elements $e$ of $\mathscr N$ are mapped to the operation $e'$
defined as follows: $e'(c)=c$, and for $x \in D \setminus \{c\}$ we set 
$e'(x) = i(e(i^{-1}(x)))$. 
All elements outside $\mathscr N$ are mapped to the constant operation $x \mapsto c$. As $\mathscr N$ is a submonoid of $\mathscr M$ and by item $(2)$, $\xi$ is a monoid homomorphism into $\mathscr O^{(1)}_D$. But $\xi$ 
is not continuous because of the first condition: if $f\in\mathscr M\setminus \mathscr N$ is contained in the closure of $\mathscr N$, then $\xi(f)$ is constant with value $c$. Thus $\xi(f)$ is not in the closure of $\xi[\mathscr N]$, as no function in $\xi[\mathscr N]$ attains the value $c$ on $D \setminus \{c\}$.
\end{proof}

Note that closed submonoids of $\mathscr O^{(1)}_D$ have certain natural submonoids, e.g., the invertible functions, or the surjective functions. Often, this yields a situation where Proposition~\ref{prop:no-automatic-cont} applies, for example in the following corollary.

\begin{cor}\label{cor:no-automatic-cont-selfembeddings-oligomorphic}
No monoid of self-embeddings of an $\omega$-categorical structure has automatic continuity.
\end{cor}
\begin{proof}
 Let $\Delta$ be an $\omega$-categorical structure on a countably infinite domain $D$, and let $\mathscr M$ be the monoid of self-embeddings of $\Delta$; we may assume that the language of $\Delta$ is countable. 
We apply Proposition~\ref{prop:no-automatic-cont} for the submonoid $\mathscr N$ of surjective functions in $\mathscr M$. We only have to show that the closure of ${\mathscr N}$ contains a non-surjective function. To prove this, observe that by the compactness theorem, $\Delta$ has a countable elementary expansion $\Delta'$ whose domain properly contains $D$.  Now $\Delta$ and $\Delta'$ are isomorphic by $\omega$-categoricity, and any isomorphism from $\Delta'$ to $\Delta$ is a non-surjective elementary self-embedding of $\Delta'$. Since $\Delta'$ has such a self-embedding, so does $\Delta$. This elementary self-embedding of $\Delta$ is contained in the closure of the automorphisms of $\Delta$, and in particular in the closure of $\mathscr N$.
\end{proof}

\subsection{Positive results: from groups to monoids}
We will focus in the following on monoids $\bf M$ with a dense subset of invertible elements. 
Our first results are a general technique
for proving automatic homeomophicity for such monoids: the basic idea is to reduce the task
to questions about certain endomorphisms
of the monoid $\bf M$. 

We then present lemmata that perform this
analysis of the endomorphisms of $\bf M$
under certain assumptions on $\bf M$;
in particular, we will assume that
the permutation group of invertible elements
of $\bf M$ has no algebraicity and 
the joint extension property.

%
The following proposition has already been outlined by Lascar in~\cite{Lascar}.

\begin{proposition}\label{Lascar}
Let $\mathbf M$ and $\mathbf M'$ be closed submonoids of $\mathbf O^{(1)}$ with dense subsets of invertibles $\mathbf G$ and $\mathbf G'$. Let $\xi\colon \mathbf G\To\mathbf G'$ be a continuous homomorphism. Then:
\begin{enumerate}
\item $\xi$ extends to a continuous homomorphism $\bar\xi\colon \mathbf M\To \mathbf M'$;
\item if $\xi$ is an isomorphism, then $\overline\xi\colon\mathbf M\To \mathbf M'$ is an isomorphism and a homeomorphism.
\end{enumerate}
\end{proposition}
\begin{proof}
 Let $\mathscr M, \mathscr M'$ be actions of $\mathbf M, \mathbf M'$ on a countably infinite set $D$, and write $\mathscr G,\mathscr G'$ for the corresponding groups of invertible functions. 
 We first show that $\xi$, as a function from $\mathscr G$ to $\mathscr G'$, is uniformly continuous with respect to the metric $d$ of Subsection~\ref{sect:topo-clones}. Let $\varepsilon>0$ be given. By continuity, there exists $\delta>0$ such that $d(\id_D,g)<\delta$ implies $d(\id_D,\xi(g))<\varepsilon$ for all $g\in \mathscr G$, where $\id_D$ denotes the identity function in $D$. Note that $d$ is a metric on $\mathscr S_D$ that is invariant under composition from the left, i.e., $d(h\circ g_1,h\circ g_2)=d(g_1,g_2)$ for all $g_1,g_2,h\in\mathscr S_D$. Hence, $d(g_1,g_2)=d(\id_D, g_1^{-1}\circ g_2)$ for all $g_1,g_2\in\mathscr S_D$. We conclude that whenever $g_1,g_2\in\mathscr G$ and $d(g_1,g_2)<\delta$, then $d(\xi(g_1),\xi(g_2))=d(\id_D,\xi(g_1)^{-1}\circ\xi(g_2))<\varepsilon$.
 
Since $\xi$ is uniformly continuous, it extends to a continuous mapping $\bar\xi\colon \mathscr M\To \mathscr M'$. Identifying the elements of ${\mathscr M}$ with equivalence classes of Cauchy sequences in ${\mathscr G}$ in the natural way, the mapping $\overline{\xi}$ sends equivalence classes of Cauchy sequences in ${\mathscr G}$ to such classes in ${\mathscr G'}$. Via this identification one easily sees that $\overline{\xi}$ is a homomorphism; this has been explictly verified in~\cite{Lascar}. This proves item~(1).

To show~(2), note that if $\xi$ is in addition an isomorphism, then it is a homeomorphism by Proposition~\ref{prop:homeo-from-cont}. In this situation, it is clear from the identification in~\cite{Lascar} that $\overline{\xi}$ is bijective, and that $(\overline{\xi})^{-1}=\overline{(\xi^{-1})}$, because Cauchy sequences in ${\mathscr G}$ correspond to Cauchy sequences in ${\mathscr G'}$ in a one-to-one manner. Hence, 
$\overline{\xi}$ is an isomorphism and a homeomorphism.
\end{proof}

\begin{lemma}\label{recmon}
 Let ${\bf M}$ be a closed submonoid of ${\bf O}^{(1)}$ whose group of invertible elements $\mathbf G$ is dense in $\mathbf M$ and has automatic homeomorphicity. Assume that the only injective endomorphism of $\bf M$ that fixes every element of $\bf G$ is the identity function $\id_{\mathbf M}$ on $\mathbf M$. Then $\bf M$ has automatic homeomorphicity.
\end{lemma}
\begin{proof}
 Let ${\bf M}'$ be a closed submonoid of ${\bf O}^{(1)}$, and let $\xi \colon \bf M \rightarrow M'$ be an isomorphism. Writing $\mathbf G'$ for the set of invertible elements of $\mathbf M'$, we have that $\xi\rest_{\bf G}$ is an isomorphism between $\mathbf G$ and $\mathbf G'$. The group $\bf G$ (and likewise $\mathbf G'$) is a closed subgroup of $\mathbf S$: when we view $\bf M$ as a closed subset of $\mathbf O^{(1)}$, and $\mathbf S$ as the subset of invertibles of $\mathbf O^{(1)}$, then $\mathbf G=\mathbf M\cap\mathbf S$. 
 As $\bf G$ has automatic homeomorphicity, we have that $\xi\rest_{\bf G}$ is a homeomorphism between $\bf G$ and $\bf G'$. Hence, by Proposition~\ref{Lascar} it extends to a mapping $\overline{\xi\rest_{\bf G}}$ from $\bf M$ to the closure $\overline{{\bf G'}}$ of $\mathbf G'$ in $\mathbf M'$; this extension is an isomorphism and homeomorphism between $\bf M$ and $\overline{{\bf G'}}$.

Let $\Phi:=\xi^{-1}\circ \overline{\xi\rest_{\bf G}}$. Then $\Phi\in \End(\bf M)$ is injective and fixed $\mathbf G$ pointwise. Thus $\Phi$ is the identity on $\mathbf M$, and consequently $\xi=\overline{\xi\rest_{\bf G}}$. In particular, $\bf G'$ is dense in $\bf M'$, and $\xi$ is a homeomorphism between $\bf M$ and $\bf M'=\overline{\mathbf G'}$. 
\end{proof}

Following \Fresse~(see~\cite{HodgesLong}), 
the \emph{age} of a relational structure $\Delta$ 
is the class of all finite structures 
that embed into $\Delta$, and denoted by $\Age(\Delta)$. 

\begin{definition}\label{def:rich}
Let $\Delta$ be a relational structure. We call a subset $U$
of the domain of $\Delta$ \emph{rich} iff for 
 every embedding $a\colon \Gamma\To \Delta$, where $\Gamma\in\Age(\Delta)$ is finite, and every $p$ in the domain of $\Gamma$ there exists an embedding $b\colon \Gamma\To\Delta$ such that $b(p) \in U$, and which agrees with $a$ on all other elements of the domain of $\Gamma$. 
We call a subset of the domain of 
$\Delta$ \emph{co-rich}
iff its complement in $\Delta$ is rich.  
\end{definition}

When $\Delta$ is a countable homogeneous
structure, then a simple back-and-forth argument shows that the structure induced by a rich subset of the domain of $\Delta$ is isomorphic to $\Delta$. 

The following definition of the concept of
\emph{no algebraicity} of permutation groups
has been given in~\cite{Oligo}. 
When the permutation group under consideration is the automorphism group of an $\omega$-categorical structure $\Delta$, then this definition coincides 
with the model-theoretic definition 
of no algebraicity for $\Delta$ (see, e.g.,~\cite{Hodges}). 

\begin{definition}
A permutation group $\mathscr G$
is said to have \emph{no algebraicity} 
iff for every finite tuple $(a_1,\dots,a_n)$ of elements of the domain of $\mathscr G$,
the set of all permutations of $\mathscr G$ that
fix each of $a_1,\dots,a_n$ fixes no other elements 
of the domain.
\end{definition}

Note that if $\mathscr G$ is a permutation group that has no algebraicity,
and $a_1,\dots,a_n$ are elements of the domain 
of $\mathscr G$, then the group of all permutations
in $\mathscr G$ that fix each of $a_1,\dots,a_n$
has only infinite orbits, except for the orbits of
$a_1,\dots,a_n$. 

\begin{lemma}\label{lem:rich}
For a countable homogeneous relational structure $\Delta$ the following are equivalent: 
\begin{itemize}
\item $\Aut(\Delta)$ has no algebraicity; 
\item $\Delta$ has a rich and co-rich subset.  
\end{itemize}
\end{lemma}

\begin{proof}
Suppose first that $\Aut(\Delta)$ 
has no algebraicity. Let $\sigma$ be the expansion of the signature $\tau$ of $\Delta$ by a new unary relation symbol $U$, and let $\mathcal C$ be the class of all finite $\sigma$-structures whose $\tau$-reduct is in 
$\Age(\Delta)$. 
Then $\mathcal C$ is a Fra\"{\i}ss\'e class.
We only indicate how to verify the amalgamation property of $\mathcal C$. Let $\Gamma_0, \Gamma_1, \Gamma_2 \in\mathcal C$, and let $s_1\colon \Gamma_0\To\Gamma_1$ and $s_2\colon \Gamma_0\To\Gamma_2$ be embeddings. By 
homogeneity of $\Delta$ there exists a 
finite substructure $\Gamma_3'$ of $\Delta$, and embeddings $t_1,t_2$ of the $\tau$-reducts of $\Gamma_1$ and $\Gamma_2$ into $\Gamma_3'$ such that $t_1\circ s_1=t_2\circ s_2$. 
It is known (see (2.15) in~\cite{Oligo}) 
that the automorphism group 
of a countable homogeneous (but not necessarily $\omega$-categorical)
structure $\Delta$ has no algebraicity if and only if
the age of $\Delta$ has \emph{strong amalgamation}. 
That is, $t_1$ and $t_2$ can be chosen such 
that $t_1[\Gamma_1] \cap t_2[\Gamma_2] = t_1[s_2[\Gamma_0]]$. Therefore, there
exists an expansion $\Gamma_3$ of $\Gamma_3'$
such that $t_1$ and $t_2$ are even embeddings
of $\Gamma_1$ and $\Gamma_2$ into $\Gamma_3$, showing the amalgamation property
for $\mathcal C$. 
We write $\Gamma$ for the \Fresse-limit for $\mathcal C$.
It can be shown by a back-and-forth argument that the $\tau$-reduct of $\Gamma$ is homogeneous and has the same age as $\Delta$. By \Fresse's theorem (see~\cite{Hodges}) it is isomorphic to $\Delta$, so let us identify the $\tau$-reduct of $\Gamma$ with $\Delta$. It is straightforward to verify that the set denoted by $U$ in $\Gamma$ is rich and co-rich with respect to $\Delta$.

For the converse implication, observe that
if $\Aut(\Delta)$ has algebraicity, 
then there is a tuple $(u,v_1,\dots,v_n)$ 
of elements of $\Delta$
such that $u$ does not appear in $\{v_1,\dots,v_n\}$
and is fixed by all automorphisms of $\Delta$
that fix each of $v_1,\dots,v_n$.
When $u \in U$ for a subset $U$ of the elements of $\Delta$, consider
the structure $S$ induced by $\{u,v_1,\dots,v_n\}$ in $\Delta$. If there were an 
embedding $b$ of $S$
into $\Delta$ such that $b(u)$ is in the complement
$U'$ of $U$ in $\Delta$,
and $b(v_i) = v_i$ for all $1 \leq i \leq n$,
then by the homogeneity of $\Delta$ 
there exists an automorphism of $\Delta$
extending $b$, contradicting the fact that $u$ is
fixed by all automorphisms of $\Delta$ which fix
each of $v_1,\dots,v_n$. Hence, $U'$ is not rich. 
\end{proof}

In the following, a \emph{partial isomorphism} 
of a structure $\Delta$ is an isomorphism between finite substructures of $\Delta$. 
We write $\Dom(a)$ for the domain 
and $\Ima(a)$ for the image of a partial isomorphism $a$.

\begin{lemma}
\label{lem:extension}
Let $\Delta$ be a countable homogeneous
relational structure such that $\Aut(\Delta)$ has no algebraicity. 
Let $f \in \overline{\Aut(\Delta)}$ have rich and co-rich image, and suppose that $a,b$ are partial isomorphisms of $\Delta$ such that 
\begin{itemize}
\item $\Dom(a) \cap \Ima(f) = f[\Ima(b)]$,
\item $\Ima(a) \cap \Ima(f) = f[\Dom(b)]$, and
\item $afb(x)=f(x)$ for all $x \in \Dom(b)$.
\end{itemize}
Then $a$ and $b$ can be extended to automorphisms $\alpha,\beta$ of $\Delta$ such that $\alpha f \beta = f$.
\end{lemma}
\begin{proof}
We construct $\alpha$ and $\beta$ by a back-and-forth argument, extending $a$ and $b$ in turns in one of the following four ways.
In each step, the three conditions on $a$ and $b$ given in the statement will be preserved.  
We write $D$ for the domain of $\Delta$. 
\begin{enumerate}
\item \emph{Extending the domain of $b$.} Let $u \in D \setminus \Dom(b)$ be arbitrary. 
Since $\Ima(a) \cap \Ima(f) = f[\Dom(b)]$ and $u \notin \Dom(b)$, we have that $f(u) \notin \Ima(a)$. 
Therefore, and since $\Delta$ is homogeneous, $a$ has an extension to a partial isomorphism $a'$ of $\Delta$ whose domain additionally contains a new element $s$ such that $a'(s) = f(u)$. 
Since $\Ima(f)$ is rich, $s$ can be chosen from $\Ima(f)$.
Let $b'$ be the extension of $b$ to the new element $u$ such that $b'(u) = f^{-1}(s)$. Then $b'$ is a partial isomorphism of $\Delta$, and $a'fb'(x)=f(x)$ for all $x \in \Dom(b')$. 
Moreover, 
\begin{align*}
\Dom(a') \cap \Ima(f) & = (\Dom(a) \cap \Ima(f)) \cup \{s\} = f[\Ima(b')] \quad \text{and} \\
\Ima(a') \cap \Ima(f) & = (\Ima(a) \cap \Ima(f)) \cup \{f(u)\} = f[\Dom(b')] \; .
\end{align*} 
\item \emph{Extending the image of $b$.} Let $v \in D \setminus \Ima(b)$ be arbitrary. 
Since $\Delta$ is homogeneous, $a$ has an extension to a partial isomorphism $a'$ of $\Delta$ with domain $\Dom(a) \cup \{f(v)\}$.
Since $\Ima(f)$ is rich, $t :=a'(f(v))$ can be chosen from $\Ima(f)$.
Let $b'$ be the extension of $b$ to the new element $f^{-1}(t)$ such that $b'(f^{-1}(t)) = v$. Then $b'$ is an isomorphism of $\Delta$, and $a'fb'(x)=f(x)$ for all $x \in \Dom(b')$. 
The other two conditions for $a'$ and $b'$ can be verified
analogously to the previous case. 

\item \emph{Extending the domain of $a$.} 
Let $s \in D \setminus \Dom(a)$. 

 Case (3.1): $s \in \Ima(f)$.
 By the homogeneity of $\Delta$ the partial isomorphism
 $a$ can be extended to a partial isomorphism $a'$ of $\Delta$ that is additionally defined on $s$; since $\Ima(f)$ is rich, we can even find an extension such that $a'(s) \in \Ima(f)$. 
Then the extension $b'$ of $b$ to $f^{-1}(a'(s))$ such that $b'(f^{-1}(a'(s))) = f^{-1}(s)$ is a partial isomorphism of $\Delta$, and $a'$ and $b'$
clearly satisfy the three conditions on $a$ and $b$ given in the statement. 

 Case (3.2): $s\not\in \Ima(f)$. 
 By the homogeneity of
 $\Delta$ the partial isomorphism $a$ can be extended to a partial isomorphism $a'$ of $\Delta$ that is additionally defined on $s$; since 
 $\Ima(f)$ is co-rich,
 we can even find an extension such that
 $a'(s) \not\in \Ima(f)$. Then $a'fb'(x)=f(x)$ for all $x \in \Dom(b)$. 
 Moreover, $\Dom(a') \cap \Ima(f) = \Dom(a) \cap \Ima(f) = f[\Ima(b)]$ and $\Ima(a') \cap \Ima(f) = \Ima(a) \cap \Ima(f) = f[\Dom(b)]$. 
 
\item \emph{Extending the image of $a$.} Let $t \in D \setminus \Ima(a)$.
 
 Case (4.1): $t\in \Ima(f)$.
 Similar to Case (3.1). 
 
 Case (4.2): $t \not\in \Ima(f)$. By the homogeneity
 of $\Delta$ the partial isomorphism $a$ can be extended to a partial isomorphism $a'$ of $\Delta$ whose domain additionally contains an element 
 $s$ such that $a'(s) = t$. 
 Since 
 $\Ima(f)$ is co-rich, 
 $s$ can be chosen 
 from $D \setminus \Ima(f)$.
 Then $a'fb(x) = f(x)$ for all $x \in \Dom(b)$. 
 Moreover, $\Dom(a') \cap \Ima(f) = \Dom(a) \cap \Ima(f) = f[\Ima(b)]$ and $\Ima(a') \cap \Ima(f) = \Ima(a) \cap \Ima(f) = f[\Dom(b)]$. 
 \end{enumerate}
 By infinite repetition of those four extension steps in turns, we arrive at automorphisms $\alpha, \beta$ of $\Delta$ such that $\alpha f \beta = f$. 
\end{proof}

\begin{lemma}\label{lem:idF}
 Let $\Delta$ be a countable
 homogeneous relational 
 structure such that $\Aut(\Delta)$ has no algebraicity. 
  Let $\mathscr F\subseteq \overline{\Aut(\Delta)}$ be the set of all self-embeddings of $\Delta$ with rich and co-rich image. Let $\xi$ be an injective endomorphism of the monoid $\overline{\Aut(\Delta)}$ which fixes $\Aut(\Delta)$ pointwise. Then $\xi$ fixes $\mathscr F$ pointwise.
\end{lemma}
\begin{proof}
Let $f \in \mathscr F$ and define $S(f) := \{(\alpha,\beta)\in \Aut(\Delta)^2 \; | \; \alpha f \beta=f\}$. 
Let $u$ and $s$ be elements of $\Delta$ 
with $s \neq f(u)$. 
We first construct a pair $(\alpha,\beta) \in S(f)$ 
 such that $\beta(u)=u$ and $\alpha(s) \neq s$. 
 
 \begin{itemize}
 \item Case 1: $s \in \Ima(f)$. Let $u'$ be the preimage of $s$ under $f$; by assumption, 
 $u' \neq u$. 
 Since $\Ima(f)$ is rich,  
 there exists an element $t \in \Ima(f)$ and
a partial isomorphism $a$ of $\Delta$ 
such that $\Dom(a) = \{f(u),s\}$, 
 $\Ima(a) = \{f(u),t\}$, 
 $a(f(u))=f(u)$, and $a(s) = t$. 
 We can additionally require that $t \in \Ima(f)$ 
 is distinct from $s$: to see this, observe that
 there exists an $s'$ distinct from $s$ such that $s$ and $s'$ lie in
 the same orbit in the group of all automorphisms of $\Delta$ that fix $f(u)$, since
 $\Aut(\Delta)$ has no algebraicity. Let $S$ be the structure induced by
 $\{f(u),s,s'\}$. Since $\Ima(f)$ is rich
 there exists an embedding $e$ of $S$ into $\Delta$
 such that $e(s) \in \Ima(f)$, and applying
 the definition of richness of $\Ima(f)$ another time,
 there also exists an embedding $e'$ of $S$ into $\Delta$
 such that $\{e'(s),e'(s')\} \subseteq \Ima(f)$. 
 Now, at least one of $e'(s)$ and $e'(s')$ is distinct from $s$. 
 
 The mapping $b$ such that 
  $\Dom(b) = \{u,f^{-1}(t)\}$, $\Ima(b) = \{u,u'\}$,
 $b(u)=u$, and $b(f^{-1}(t)))=u'$ is a partial isomorphism of $\Delta$. Moreover, $a$ and $b$ satisfy the conditions of Lemma~\ref{lem:extension}: 
we have $\Dom(a) \cap \Ima(f) = \{f(u),s\} = \{f(u),f(u')\} = f(\{\Ima(b)\})$, 
$\Ima(a) \cap \Ima(f) = \{f(u),t\} = \{f(u),f(f^{-1}(t))\} = f(\Dom(b))$, and $afb(x) = f(x)$ for $x \in \{u,f^{-1}(t)\}$. 
\item Case 2: $s \notin \Ima(f)$. Let $b$
be the partial isomorphism of $\Delta$
such that $\Dom(b) = \{u\}$, $\Ima(b) = \{u\}$. 
Since $\Ima(f)$ is co-rich, there is an element
$t$ of $\Delta$ outside of $\Ima(f)$
and a partial isomorphism $a$ of $\Delta$
such that $\Dom(a) = \{f(u),s\}$, $\Ima(a) = \{f(u),t\}$, $a(f(u))=f(u)$, and $a(s) = t$. 
Using the fact that $\Ima(f)$ is co-rich, and arguing as in the previous item, we can additionally assume that $t$ is distinct from $s$.

Then $a$ and $b$
satisfy the conditions of Lemma~\ref{lem:extension}:
$\Dom(a) \cap \Ima(f) = \{f(u)\} = f(\{\Ima(b)\})$, 
$\Ima(a) \cap \Ima(f) = \{f(u)\} = f(\Dom(b))$, and $afb(u) = f(u)$. 
\end{itemize}
In both cases, by Lemma~\ref{lem:extension},
 there are 
 $(\alpha,\beta) \in S(f)$ 
 such that $\beta(u)=u$ 
 and $\alpha(s) = t \neq s$.
 
 We can now describe how the element $f(u)$ 
  can be recovered from $S(f)$, 
  namely as 
  \begin{align}
  \{f(u)\}=\bigcap\limits_{(\alpha, \beta)\in S(f)\atop \beta(u)=u} \{s \; | \; \alpha(s) = s\} \; . \label{eq:fix}
  \end{align}
  Thus, $\xi(f)=f$. 
  
 For the inclusion ``$\subseteq$'' in Equation $(\ref{eq:fix})$,
 note that
  the conditions $\beta(u)=u$ 
  and $(\alpha, \beta)\in S(f)$ 
  imply that $\alpha(f(u))=f(u)$. Hence, $f(u)$ belongs to the right-hand-side of Equation~$(\ref{eq:fix})$.
	For the inclusion ``$\supseteq$'', let $s$ be any element of $\Delta$ distinct from $f(u)$. Then we have seen above that there
	exists $(\alpha,\beta) \in S(f)$ such that
	$\beta(u)=u$ and $\alpha(s) \neq s$. 
	Hence, $s$ does not belong to the right-hand-side of Equation $(\ref{eq:fix})$. 
\end{proof}

\begin{definition}
We say that a structure $\Delta$ has the
\emph{joint extension property} 
iff for all partial isomorphisms 
$a_1,a_2$ of $\Delta$ with $\Dom(a_1) = \Dom(a_2)$ and $a_1^{-1}(x) = a_2^{-1}(x)$ for all $x \in \Ima(a_1) \cap \Ima(a_2)$,
and for every element $u$ of $\Delta$
outside of $\Dom(a_1) = \Dom(a_2)$ there exist
extensions $a_1'$ of $a_1$ and $a_2'$ of $a_2$
such that $a_1'(u) = a_2'(u)$. 
\end{definition}

Examples of structures with the joint extension property are the random graph, 
the random tournament, and the random digraph. 
Examples of structures without the joint extension property are $({\mathbb Q};<)$ and for $n \geq 3$ the countable universal homogeneous $K_n$-free graph $G$. To see this for $n = 3$, let $v$ be a vertex of $G$, and let $a_1,a_2$ be maps with
domain $\{v\}$ such that $a_1(v)$ is adjacent to $a_2(v)$. Let $u$ be a vertex adjacent to $v$
in $G$. 
Then in any extension of $a_1'$ of $a_1$
and $a_2'$ of $a_2$ with domain $\Dom(a_1')=\Dom(a_2') = \{u,v\}$ such that $a_1'(u)=a_2'(u)$ we must
have that $a_1'(u)$ is not adjacent to $a_1'(v)$ or $a_2'(u)$ is not adjacent to $a_2'(v)$, 
since otherwise the three vertices
$a_1'(v),a_1'(u)=a_2'(u),a_2'(v)$ would form a triangle. Hence, $a_1$ and $a_2$ cannot be embeddings, showing that the joint extension property fails. 

\begin{lemma}\label{Tu}
 Let $\Delta$ be a countable homogeneous relational structure 
 and with the joint extension property such that $\Aut(\Delta)$ has no algebraicity.
 Let $f\in \overline{\Aut(\Delta)}$, and suppose that $a_1,a_2,b$ are 
 partial isomorphisms 
 of $\Delta$ such that 
 \begin{itemize}
 \item $f[\Ima(b)] \subseteq \Dom(a_1) = \Dom(a_2)$,
 \item $a_1^{-1}(x) = a_2^{-1}(x)$ for all $x \in \Ima(a_1) \cap \Ima(a_2)$, and
\item $a_1 f b(x) = a_2 f b(x)$ for all $x \in \Dom(b)$. 
\end{itemize}
Then $a_1,a_2,b$ extend to  
$\alpha_1,\alpha_2,\beta\in \overline{\Aut(\Delta)}$ with rich and co-rich image such that $\alpha_1 f \beta = \alpha_2 f \beta$. 
\end{lemma}
\begin{proof}
We construct $\alpha_1$, $\alpha_2$, and $\beta$ 
by a back-and-forth argument, stepwise extending
$a_1$, $a_2$, and $b$.
In our construction, we also maintain finite subsets $A_1$, $A_2$, $B$ of $D$
such that at each stage during the construction, $\Ima(a_1) \cap A_1 = \emptyset$,
$\Ima(a_2) \cap A_2 = \emptyset$, and
$\Ima(b_1) \cap B = \emptyset$. 
Initially, we set $A_1=A_2=B= \emptyset$. 
In each step, we either extend $a_1$, $a_2$, and $b$ to partial isomorphisms of $\Delta$ 
such that the three conditions from the statement on $a_1,a_2,b$ remain valid, or we add elements to
$A_1,A_2$, and $B$ to make sure that the images of $\alpha_1$, $\alpha_2$ and $\beta$ will be
co-rich. 

\begin{enumerate}
\item \emph{Extending the domain of $b$.} 
Let $u  \in D \setminus \Dom(b)$ be arbitrary. 
Since $\Aut(\Delta)$ is without algebraicity, there exists
an element $v \in D \setminus (B \cup f^{-1}[\Dom(a_1)])$ and an extension $b'$ of $b$ to a partial isomorphism 
of $\Delta$ such that $b'(u) = v$. 
By the joint extension property of $\Delta$,
there are partial isomorphisms $a_1'$ extending $a_1$
and $a_2'$ extending $a_2$ that are additionally defined
on $f(v)$ and $a_1'(f(v))=a_2'(f(v))$.
Since $\Aut(\Delta)$ has no algebraicity
we can assume that $a_1'(f(v)) = a_1'(f(v)) \notin A_1 \cup A_2$.
\item \emph{Extending the domain of $a_1$ and $a_2$.} 
Let $s \in D \setminus \Dom(a_1)$.
Since $\Delta$ is homogeneous and $\Aut(\Delta)$ is without algebraicity, there is an extension of $a_1$ to a partial isomorphisms 
$a_1'$ of $\Delta$ which 
is additionally defined on $s$ 
such that $a_1'(s) \notin (A_1 \cup \Ima(a_2))$.
Similarly, there is an extension of $a_2$ to a partial isomorphisms 
$a_2'$ of $\Delta$ which 
is additionally defined on $s$ 
such that $a_2'(s) \notin (A_2 \cup \Ima(a'_1))$.
\item \emph{Extending $B$.}
Pick $d \in D \setminus (B \cup \Ima(b))$,
and add $d$ to $B$. 
\item \emph{Extending $A_1$ and $A_2$.}
Pick $d \in D \setminus (A_1 \cup \Ima(a_1))$,
and add $d$ to $A_1$. We extend $A_2$ similarly. 
\item \emph{Enriching the image of $b$.} 
Let $S \in \Age(\Delta)$, let $p \in S$, and
let $e$ be an embedding of $S$ into $\Delta$ such that $e(x) \in B \cup \Ima(b)$ for all $x \in S \setminus \{p\}$.
Since $\Delta$ is homogeneous
and $\Aut(\Delta)$ is without algebraicity, there is an element $q \in D \setminus (\Ima(b) \cup B)$ and an embedding
$e'$ of $S$ into $\Delta$ such that
$e'(p) = q$ and $e'(x) = e(x)$ for all 
$x \in S \setminus \{p\}$. 
By the homogeneity of $\Delta$, there is
$u \in D$ and an extension of $b$ to a partial isomorphism $b'$ of $\Delta$ with $\Dom(b') = \Dom(b) \cup \{u\}$ such that $b'(u) = q$.
Since $\Aut(\Delta)$ has no algebraicity and the joint extension property, 
there exist extensions $a_1'$ and $a_2'$ 
of $a_1$ and $a_2$
that are additionally defined on $f(q)$ such that
$a_1'(f(q)) = a_2'(f(q)) \notin A_1 \cup A_2$. 
%
\item \emph{Enriching the image of $a_1$ and $a_2$.} 
Let $S \in \Age(\Delta)$, and $p \in S$, and
$e$ be an embedding of $S$ into $\Delta$ such that
$e(x) \in \Ima(a_1) \cup A_1$ for all $x \in S \setminus \{p\}$.
Since $\Delta$ is homogeneous
and $\Aut(\Delta)$ has no algebraicity, there is an element $q \in D \setminus (\Ima(a_1) \cup A_1)$ 
and an embedding
$e'$ of $S$ into $\Delta$ such that
$e'(p) = q$ and $e'(x) = e(x)$ for all 
$x \in S \setminus \{p\}$. 
By the homogeneity of $\Delta$
there is an element $s \in D \setminus \Dom(a_1)$ and an extension $a_1'$ of $a_1$ to
a partial isomorphism of $\Delta$ that is 
additionally defined on $s$ such that $a_1'(s) = q$.
Since $\Delta$ is homogeneous and $\Aut(\Delta)$ is without algebraicity, $a_2$ has an extension $a_2'$ to
a partial isomorphism of $\Delta$ that is additionally defined on $s$ such that $a_1'(s) \notin A_1$. We extend $a_2$ similarly. 
\item \emph{Enriching $B$.} 
Let $S \in \Age(\Delta)$, and $p \in S$, and
$e$ an embedding of $S$ into $\Delta$ such that
$e(x) \in \Ima(b) \cup B$ for all $x \in S \setminus \{p\}$.
Since $\Delta$ is homogeneous
and $\Aut(\Delta)$ without algebraicity, there is an element $q \in D \setminus (\Ima(b) \cup B)$ and an embedding
$e'$ of $S$ into $\Delta$ such that
$e'(p) = q$ and $e'(x) = e(x)$ for all 
$x \in S \setminus \{p\}$. Add $q$ to $B$.
%
\item \emph{Enriching $A_1$ and $A_2$.} 
Let $S \in \Age(\Delta)$, let $p \in S$, and
let $e$ an embedding of $S$ into $\Delta$ such that
$e(x) \in \Ima(a_1) \cup A_1$ for all $x \in S \setminus \{p\}$.
Since $\Delta$ is homogeneous
and $\Aut(\Delta)$ without algebraicity, there is an element $q \in D \setminus (\Ima(a_1) \cup A_1)$ 
and an embedding
$e'$ of $S$ into $\Delta$ such that
$e'(p) = q$ and $e'(x) = e(x)$ for all 
$x \in S \setminus \{p\}$. 
Add $q$ to $A_1$. 
Similarly, we extend $A_2$.
\end{enumerate}
We perform these steps in turns so that
each step is performed infinitely often. 
Because of $(1)$ and $(2)$ we can make sure that we arrive at self-embeddings
$\alpha_1,\alpha_2,\beta$ of $\Delta$ such that
$\alpha_1 f \beta = \alpha_2 f \beta$. 
Step $(3)$ makes sure
that the union over the sets $B$ from the construction equals the complement of $\Ima(\beta)$. Similarly, step $(4)$ makes sure that
the union of the $A_1$ and over the $A_2$ 
yield the complement of $\Ima(\alpha_1)$
and $\Ima(\alpha_2)$, respectively. 
Repeating step $(5)$ infinitely often 
for all possible $S \in \Age(\Delta)$ ensures that $\Ima(\beta)$ is rich. 
Repeating step $(6)$ infinitely often for all possible $S \in \Age(\Delta)$ ensures that $\Ima(\alpha_1)$ and $\Ima(\alpha_2)$ are rich. 
Repeating step $(7)$ infinitely often for all possible $S \in \Age(\Delta)$ ensures that the complement 
of $\Ima(\beta)$ is rich. 
Repeating step $(8)$ infinitely often for all possible $S \in \Age(\Delta)$
makes sure that the complement of $\Ima(\alpha_1)$ and the complement of $\Ima(\alpha_2)$ are rich. Since $\alpha_1,\alpha_2$, and $\beta$ are embeddings, and every set that contains a rich set is rich, 
the images of $\alpha_1$, $\alpha_2$, and $\beta$
are rich, too. 
\end{proof}

\begin{lemma}\label{lem:idM}
 Let $\Delta$ be a countable homogeneous relational structure with the joint extension property such that $\Aut(\Delta)$ has no algebraicity. 
 Let $\xi$ be an injective endomorphism of the monoid $\overline{\Aut(\Delta)}$ which fixes $\Aut(\Delta)$ pointwise. 
 Then $\xi$ is the identity.
\end{lemma}
\begin{proof}
  We write $\mathscr F$ for the set of self-embeddings of $\Delta$ with rich and co-rich image.
  For any element $f$ of $\overline{\Aut(\Delta)}$ 
  define $T(f) := \{(\alpha_1, \alpha_2, \beta) \in {\mathscr F}^3 \; | \; \alpha_1 f \beta = \alpha_2 f \beta\}$.
 Let $u,s$ be elements of $\Delta$ with $s \neq f(u)$. Our goal is to construct a triple $(\alpha_1, \alpha_2, \beta) \in T(f)$ such that $\beta(u)=u$ and $\alpha_1(s) \neq \alpha_2(s)$.
Since $\Aut(\Delta)$ has no algebraicity,
there are two distinct elements $t_1,t_2$
and partial isomorphisms $a_1,a_2$
of $\Delta$ with $\Dom(a_1) = \Dom(a_2) = \{f(u),s\}$, $a_1(f(u))=a_2(f(u))=f(u)$,
$a_1(s)=t_1$, and $a_2(s) = t_2$.
Let $b$ be the partial isomorphism such that $\Dom(b) = \Ima(b) = \{u\}$. 
Then the conditions of Lemma~\ref{Tu} are satisfied, and we obtain the desired triple $(\alpha_1, \alpha_2, \beta)$.
The conditions $\beta(u)=u$ and $(\alpha_1, \alpha_2, \beta) \in T(f)$ imply that $\alpha_1(f(u))=\alpha_2(f(u))$. Hence, the element $f(u)$ can be recovered from $T(f)$, namely $$\{f(u)\}=\big \{s \in \Delta \; | \; \text{ for all } (\alpha_1, \alpha_2, \beta) \in T(f) \text{ with } \beta(u)=u
  \text{ it holds that } \alpha_1(s)=\alpha_2(s) \big \} \; .$$ 
Thus $\xi(f)=f$. 
\end{proof}

\begin{theorem}\label{thm:mainmon}
 Let $\Delta$ be a countable homogeneous 
 relational 
 structure such that $\Aut(\Delta)$ has no algebraicity and with the joint extension property such that 
 $\Aut(\Delta)$ has automatic homeomorphicity. 
 Then the monoid $\overline{\Aut(\Delta)}$ of self-embeddings of $\Delta$ has automatic homeomorphicity.
 \end{theorem}
 
\begin{proof}
By Lemma~\ref{lem:idM}, every injective endomorphism $\xi$ of the monoid $\overline{\Aut(\Delta)}$ with $\xi\upharpoonright_{\Aut(\Delta)}=\id_{\Aut(\Delta)}$ is the identity. Hence, the conditions of Lemma~\ref{recmon} hold, and thus $\overline{\Aut(\Delta)}$ has automatic homeomorphicity.
\end{proof}

\begin{corollary}\label{cor:mon}
The self-embedding monoids of the following structures have automatic homeomorphicity: 
the structure without structure, $({\mathbb N};=)$;
the random tournament;
the random graph;
the random directed graph;
the random $k$-uniform hypergraph for $k \geq 2$.
 \end{corollary}
\begin{proof}
It is well-known that all structures that appear in the statement are homogeneous structures and
have automorphism groups without algebraicity. It is easy to verify that all these structures have the joint extension property. 
References for the proofs of automatic continuity of the respective groups have been given in Section~\ref{sect:groups}; automatic homeomorphicity of the groups follows from Proposition~\ref{prop:homeo-from-cont}. 
Finally, automatic homeomorphicity of the given monoids follows from Theorem~\ref{thm:mainmon}.
\end{proof}

\section{Topological Clones}
\label{sect:clones}

\subsection{Birkhoff's theorem and continuity}\label{subsect:birkhoffcont}

\begin{defn}
If $\mathscr C$ is a function clone acting on a set $C$, then we write $(C;\mathscr C)$ for any algebra on $C$ whose fundamental operations are precisely the operations of $\mathscr C$. Note that there are many such algebras, depending on the indexing of the functions in $\mathscr C$; we emphasize that we also allow multiple appearances of the same function in $\mathscr C$ in an indexing.
\end{defn}

When $\mathcal C$ is a class of algebras with common signature $\tau$, then $\PPP(\mathcal C)$ denotes the class of all products of algebras from $\mathcal C$, 
$\SSS(\mathcal C)$ denotes the class of all subalgebras of algebras from $\mathcal C$, and $\HHH(\mathcal C)$ denotes the class of all 
homomorphic images of algebras from $\mathcal C$. The following classical theorem from universal algebra gives us representations of all (not necessarily faithful) actions of the abstract clone of a given function clone $\mathscr C$ by means of these operators.

\begin{thm}[Birkhoff~\cite{Bir-On-the-structure}]\label{thm:birkhoff}
Let $\mathscr C, \mathscr D$ be function clones acting on sets $C$, $D$ respectively. Then there exists a surjective homomorphism $\xi\colon \mathscr C\To \mathscr D$ if and only if $(D;\mathscr D)\in \HSP(C;\mathscr C)$ for some indexing of those algebras.
\end{thm}


\begin{thm}\label{prop:Oreconstruction}
Any function clone $\mathscr C$ with domain $C$ which contains $\mathscr O^{(1)}_C$ has automatic continuity.
\end{thm}
\begin{proof}
Let $D$ be a countably infinite set, and let $\xi\colon \mathscr C\To\mathscr D$ be a homomorphism onto a (not necessarily closed) subclone $\mathscr D$ of $\mathscr O_D$. Then $(D;\mathscr D)\in \HSP(C;\mathscr C)$ by Theorem~\ref{thm:birkhoff}. In other words, there is a subset $S$ of some power $C^I$ and an equivalence relation $\sim$ on $S$ such that both $S$ and $\sim$ are invariant under the componentwise action of $\mathscr C$ on $C^I$ and such that the algebra $(S;\mathscr C)/_\sim$ is isomorphic to  $(D;\mathscr D)$. The clone isomorphism $\xi$ is then obtained by sending every function $f\in \mathscr C$ to the corresponding function (with the same name) in $(S;\mathscr C)/_\sim$. 

In the following, we view tuples in $C^I$ as functions from $I$ to $C$, and in particular are going to speak of the \emph{range} and the \emph{kernel} of a tuple. Note first that since $S$ is invariant under $\mathscr O^{(1)}_C$, it contains with every tuple $a\in C^I$ all tuples whose  kernel is at least as coarse as the kernel of $a$; in other words, whether or not a tuple is in $S$ only depends on its kernel, and the kernels of tuples in $S$ are upward closed in the lattice of equivalence relations on $I$ (where the order is containment). Now let $F\subseteq S$ consist of those tuples which have finitely many values  (i.e., whose kernel has finitely many classes). Then $(F;\mathscr C)$ is a subalgebra of $(S;\mathscr C)$, as finitary operations on finite domains have finite range. We will now show that every tuple in $S$ is $\sim$-equivalent to a tuple in $F$, i.e., one with finite range. We then have that $(F;\mathscr C)/_\sim$ and $(S;\mathscr C)/_\sim$ are isomorphic. But the mapping which sends every function in $\mathscr C$ to its corresponding function in $(F;\mathscr C)/_\sim$ is continuous: to see this, note first that the mapping which sends every function in $\mathscr C$ to its corresponding function in $(F;\mathscr C)$ is continuous, since all tuples in $F$ take only finitely many values. Moreover, the mapping which sends every function in $(F;\mathscr C)$ to the corresponding function acting on the classes of $\sim$ is continuous in general, and so our claim follows by composing the two mappings. It might be worth noting that in general, among $\HHH, \SSS, \PPP$, the only operator which might act discontinuously on a function clone is the operator $\PPP$; however, it is easy to see that it does act continuously if the product is finite, or is restricted to finite range tuples.

So given $t\in S$, we show that it is $\sim$-equivalent to a tuple in $F$. We may assume that $t$ has infinite range, for otherwise there is nothing to show; in particular, $C$ is infinite. Next observe that if there exists $t'\in S$ with the same kernel as $t$ which is $\sim$-equivalent to a tuple in $F$, then $t$ is equivalent to a tuple in $F$ as well: for in that situation, there exists $f\in\mathscr O_C^{(1)}$ sending $t'$ to $t$. Hence, if $c\in F$ is so that $t'\sim c$, then $t\sim f({c})$, and $f({c})\in F$. Now there is a continuum of tuples with the same kernel as $t$ and such that the ranges of any two tuples of the continuum have finite intersection (the ranges form an \emph{almost disjoint family}, see~\cite{Jec-Set-theory}). Because $D$ is countable, $\sim$ has only countably many classes, and thus there exist $t', t''$ in the continuum with $t'\sim t''$. Pick any function $f\in\mathscr O_C^{(1)}$ which is injective on the range of $t'$ and takes only finitely many values on the range of $t''$. We then have $f(t')\sim f(t'')$, and hence $f(t')$ is equivalent to a tuple in $F$. But $f(t')$ has the same kernel as $t$, and hence also $t$ is equivalent to a tuple in $F$. 
\end{proof}

Recall from Section~\ref{sect:groups} that $\mathbf S$ has automatic continuity; we obtain the analogous statement for $\mathbf O$ and $\mathbf O^{(1)}$ as a corollary of Theorem~\ref{prop:Oreconstruction}.

\begin{cor}\label{cor:Oreconstruction}
$\mathbf O$ and $\mathbf O^{(1)}$ have automatic continuity.
\end{cor}

\subsection{Constants and openness}\label{subsect:constantsopenness}

\begin{prop}\label{prop:constantsopen}
Let $\mathscr C$ be a closed function clone with domain $C$ which contains all constant functions on $C$, and let $\xi\colon \mathscr C\To \mathscr D$ be an isomorphism onto a function clone $\mathscr D$. Then the image of any open subset of $\mathscr C$ under $\xi$ is open in $\mathscr D$.
\end{prop}
\begin{proof}
For $n\geq 1$ and $a_1,\ldots,a_n,b\in C$, let $U$ be the basic clopen set of all $n$-ary functions $f\in\mathscr C$ with $f(a_1,\ldots,a_n)=b$. Then writing $g_a$ for the unary constant function on $C$ with value $a$ for all $a\in C$, we have that $U$ consists precisely of those $n$-ary functions $f\in\mathscr C$ for which the equation $g_b(x)=f(g_{a_1}(x),\ldots,g_{a_n}(x))$ holds. Because $\xi$ is an isomorphism, $\xi[U]$ consists precisely of those $n$-ary functions $f$ in $\mathscr D$ for which the equation $\xi(g_b)(x)=f(\xi(g_{a_1})(x),\ldots,\xi(g_{a_n})(x))$ holds. Hence, $\xi[U]$ is a closed subset of $\mathscr D$ (in fact, since $\xi$ is a clone homomorphism, the functions $\xi(g_{a_k})$ are constant, so that one even sees right away that $\xi[U]$ is clopen).
\end{proof}

\begin{cor}\label{cor:aboveO1ah}
Any closed function clone $\mathscr C$ with domain $C$ which contains $\mathscr O^{(1)}_C$ has automatic homeomorphicity. In particular, $\mathbf O$ and $\mathbf O^{(1)}$ have automatic homeomorphicity.
\end{cor}
\begin{proof}
This is a direct consequence of Theorem~\ref{prop:Oreconstruction} and Proposition~\ref{prop:constantsopen}.
\end{proof}

\ignore{
\begin{prop}\label{prop:Unarypredicatesreconstruction}
Let $D$ be disjoint union of countable sets $A$ and $B$. Then $\Pol(D;A,B)$ has automatic homeomorphicity.
\end{prop}
\begin{proof}
Set $\mathscr C:=\Pol(D;A,B)$, and let $\xi\colon \mathscr C\To\mathscr D$ be an isomorphism onto a closed subclone $\mathscr D$ of $\mathscr O_D$. Then $(D;\mathscr D)\in \HSP(D;\mathscr C)$ by Theorem~\ref{thm:birkhoff}. Denote the algebra obtained by restricting the functions of $\mathscr C$ to $A$ by $(A;\mathscr C)$. It is easy to see that all algebras in $\HSP(D;\mathscr C)$ are, up to isomorphism, products of a power of $(D;\mathscr C)$ with a power of $(A;\mathscr C)$ with a power of a two-element algebra. For cardinality reasons, all these powers must be finite. Because $\mathscr C$ is oligomorphic, $\xi$  is continuous by Theorem~\ref{thm:topobirkhoff}.
\end{proof}
}
\subsection{Transitivity and openness}\label{subsect:transitivityopenness}

\begin{defn}\label{defn:transitive}
We call a function clone \emph{transitive} iff the permutation group of its invertible unary functions acts transitively on the domain of the clone. We call a topological clone transitive iff it is the topological clone of a transitive function clone. 
\end{defn}

\begin{remark}
If $\mathscr G$ is a closed oligomorphic permutation group, then the corresponding topological group $\mathbf G$ has a faithful transitive action whose image is closed. 
Indeed, if $\mathscr G$ has $k$ orbits $X_1, \ldots, X_k$, then consider the coordinate-wise action of $\mathbf G$ on $X=X_1\times \cdots \times X_k$. Pick $x_1\in X_1, \ldots, x_k\in X_k$, and let $Y$ be the orbit of $(x_1, \ldots, x_k)\in X$. 
The action of  $\mathbf G$ on $Y$ is transitive and has a closed image. 
To show faithfulness, observe that $\alpha\in \mathbf G$ is in the pointwise stabiliser of $Y$ if and only if $\alpha$ stabilises every element in the action of  $\mathscr G$ on $X_1\cup \cdots \cup X_k$, and thus $\alpha$ is the identity element. 
\end{remark}

The following example shows that contrary to closed oligomorphic subgroups of $\mathbf S$, closed oligomorphic subclones of $\mathbf O$ need not be transitive.

\begin{prop}\label{prop:notransitiveaction}
There exists an oligomorphic closed subclone of $\mathbf O$ which is not transitive.
\end{prop}
\begin{proof}
Let $D$ be the disjoint union of countable sets $A$ and $B$, and set $\mathscr C:=\Pol(D;A,B)$. Clearly, $\mathscr C$ is oligomorphic; let $\mathbf C$ be its topological clone. Suppose that $\mathbf C$ is also the topological clone of a transitive function clone $\mathscr C'$ with domain $D'$. Then $(D';\mathscr C')\in \HSP(D;\mathscr C)$ by Theorem~\ref{thm:birkhoff}. Let $I$ be a set, $S$ be a subuniverse of $(D^I;\mathscr C)$, and $\sim$ be an equivalence relation on $S$ which is invariant under $\mathscr C$ such that $(D';\mathscr C')$ is isomorphic to $(S;\mathscr C)/_\sim$. 
Consider two arbitrary equivalence classes $P,Q\subseteq S$ of $\sim$. By the transitivity of $\mathscr C'$, there exists an invertible $\alpha\in\mathscr C^{(1)}$ such that $\alpha({P})=Q$ in the interpretation of $\alpha$ in the algebra $(S;\mathscr C)/_\sim$. Now if $t$ is an arbitrary tuple in $P$, then $\alpha({P})=Q$ equals the $\sim$-class of $\alpha(t)\in S$. The tuple $\alpha(t)$ takes values in $A$ precisely when $t$ does. Therefore, if $f,g\in \mathscr C$ are binary functions which agree on $A^2$ and $B^2$, then $f(P,Q)=g(P,Q)$ in their action in $(S;\mathscr C)/_\sim$, since $f(P,Q)$ and $g(P,Q)$ is the $\sim$-class of the tuple $f(t,\alpha(t))=g(t,\alpha(t))$. Hence, $f$ and $g$ are equal in their interpretation in $(S;\mathscr C)/_\sim$, even if they differ on $A\mult B$. This contradicts the fact that the mapping which sends every $f\in\mathscr C$ to the function with the same name in $(S;\mathscr C)/_\sim$ is an isomorphism.
\end{proof}

The proof of Theorem~\ref{prop:Oreconstruction} can basically be copied to obtain that the function clone from the proof of Proposition~\ref{prop:notransitiveaction} without transitive action has automatic continuity.

\begin{prop}
Let $D$ be the disjoint union of countable sets $A$ and $B$. Then $\Pol(D;A)$ and $\Pol(D;A,B)$ have automatic continuity.
\end{prop}

We will now see how transitivity of a topological clone helps to lift openness of isomorphisms from this clone from the unary part to higher arities.

\begin{prop}\label{prop:trans-open}
Let $\mathbf C$ be a transitive topological clone, 
and let $\xi$ be an isomorphism from $\mathbf C$ onto 
a topological clone $\mathbf C'$. If the restriction of $\xi$ to $\mathbf C^{(1)}$
is open, then so is $\xi$. 
\end{prop}
\begin{proof}
For any clone $\mathfrak W$ and $g_1,\ldots,g_k\in \mathfrak W^{(1)}$, let $p_{g_1,\ldots,g_k}$ be the mapping $p_{g_1,\ldots,g_k}\colon \mathfrak W^{(k)}\To  \mathfrak W^{(1)}$ defined by 
$$p_{g_1,\ldots,g_k}(f(x_1,\ldots,x_k)):=f(g_1(x),\ldots,g_k(x)).$$
Let $\mathscr C$ be a transitive function clone acting on a set $D$ such that $\mathbf C$ is the topological clone of $\mathscr C$. Now let $k\geq 1$ and $a_0,\ldots,a_k\in D$ be given, and let $U = \{ f \in \mathscr C^{(k)} \; | \; f(a_1,\dots,a_k) = a_0\}$. The family of such subsets $U\subseteq \mathscr C$ form a subbasis of the topology on $\mathscr C$. Since $\xi$ is injective, it is enough to show that $\xi[U]$ is open. Let $\mathscr G$ be the permutation group of invertibles of $\mathscr C^{(1)}$. Since $\mathscr G$ acts transitively,  there are $\alpha_1,\dots,\alpha_k \in \mathscr G$ and 
 $b \in D$ 
such that $\alpha_i(b) = a_i$ for all $1\leq i \leq k$. Set $U':=\{g\in \mathscr C^{(1)} \; | \; g(b) = a_0\}$. We claim that
\begin{align*}
\xi[U]&=\xi[p_{\alpha_1,\ldots,\alpha_k}^{-1}[U']]=p_{\xi(\alpha_1),\ldots,\xi(\alpha_k)}^{-1}[\xi[U']].
\end{align*}
The first equation is clear since $U$ is the preimage of $U'$ under $p_{\alpha_1,\ldots,\alpha_k}$. The second equation is a direct consequence of the fact that $\xi$ is an isomorphism. 

We thus have that $\xi[U]$ is open since $\xi[U']$ is open and $p_{\xi(\alpha_1),\ldots,\xi(\alpha_k)}$ is continuous, proving the proposition.
\end{proof}

\begin{cor}\label{for:aboveSinftyopen}
Let $\mathscr C$ be a closed function clone with countably infinite domain $D$ which contains the group $\mathscr S_D$ of all permutations on $D$, and let $\xi\colon \mathscr C\To \mathscr D$ be an isomorphism onto a function clone $\mathscr D$. Then the image of any open subset of $\mathscr C$ under $\xi$ is open in $\mathscr D$.
\end{cor}
\begin{proof}
It is known that either $\mathscr C$ contains all constant functions, or its unary part consists precisely of all injections (i.e., the closure of $\mathscr S_D$ in $\mathscr O^{(1)}$)~\cite{BodChenPinsker}. In the first case the claim follows from Corollary~\ref{cor:aboveO1ah}, in the latter case from Corollary~\ref{cor:mon} and Proposition~\ref{prop:trans-open}.
\end{proof}

Proposition~\ref{prop:trans-open} can also be used in the other direction for showing continuity. Although the proof is dual, we include it for the convenience of the reader.

\begin{prop}\label{prop:cont-trans}
Let $\mathbf C$ be a topological clone,
and let $\xi$ be a homomorphism from $\mathbf C$ onto 
a transitive topological clone $\mathbf C'$. If the restriction of $\xi$ to $\mathbf C^{(1)}$
is continuous, then $\xi$ is continuous as well.  
\end{prop}
\begin{proof}
Let $D$ be the set on which $\mathbf C'$ acts transitively as a function clone $\mathscr C'$. Now let $k\geq 1$ and $a_0,\ldots,a_k\in D$ be given, and let $U = \{ f \in \mathscr C'^{(k)} \; | \; f(a_1,\dots,a_k) = a_0\}$. As in the proof of Proposition~\ref{prop:trans-open} it suffices to show that $\xi^{-1}[U]$ is open. Let $\mathscr G'$ be the permutation group of invertibles of $\mathscr C'^{(1)}$. Since $\mathscr G'$ acts transitively,  there are $\alpha_1,\dots,\alpha_k \in \mathscr G'$ and 
 $b \in D$ 
such that $\alpha_i(b) = a_i$ for all $1\leq i \leq k$. Set $U':=\{g\in \mathscr {C'}^{(1)} \; | \; g(b) = a_0\}$. Using the fact that $\xi$ is onto, pick $\beta_i\in\mathscr C$ such that $\xi(\beta_i)=\alpha_i$, for all $1\leq i\leq k$. We claim that
\begin{align*}
\xi^{-1}[U]&=\xi^{-1}[p_{\alpha_1,\ldots,\alpha_k}^{-1}[U']]=p_{\beta_1,\ldots,\beta_k}^{-1}[\xi^{-1}[U']].
\end{align*}
The first equation is clear since $U$ is the preimage of $U'$ under $p_{\alpha_1,\ldots,\alpha_k}$. To see the second equation, let $f\in \xi^{-1}[p_{\alpha_1,\ldots,\alpha_k}^{-1}[U']]$. Then $\xi(f(\beta_1,\ldots,\beta_k))=\xi(f)(\alpha_1,\ldots,\alpha_k)\in U'$, so that indeed $f\in p_{\beta_1,\ldots,\beta_k}^{-1}[\xi^{-1}[U']]$. The implications can be reversed, proving the other inclusion.

We thus have that $\xi^{-1}[U]$ is open since $\xi^{-1}[U']$ is open and $p_{\beta_1,\ldots,\beta_k}$ is continuous, proving the proposition.
\end{proof}

\subsection{Gates and continuity}\label{subsect:gatescontinuity}

\subsubsection{Simple gate coverings}

The following concept will be useful when we wish to prove continuity of a clone homomorphism knowing that it is continuous on unary functions.

\begin{defn}
A \emph{gate covering} of a topological clone $\bf C$ consists of
\begin{itemize}
\item an open covering $\mathcal U$ of $\bf C$ and
\item a function $f_U\in U$ for each $U\in {\mathcal U}$
\end{itemize}
such that for all $U\in\mathcal U$ and all Cauchy sequences $(g^j)_{j\in\omega}$ of functions in $U$ of the same arity $n\geq 1$ there exist unary functions $\alpha^j,\beta_i^j\in \bf C$, where $j\in\omega$ and $1\leq i\leq n$, such that
\begin{itemize}
\item $g^j(x_1,\ldots,x_n)=\alpha^j(f_U(\beta^j_1(x_1),\ldots,\beta^j_n(x_n)))$;
\item $(\alpha^j)_{j\in\omega}$, $(\beta_i^j)_{j\in\omega}$ are Cauchy for all $1\leq i\leq n$.
\end{itemize}
The functions $f_U$ are called the \emph{gates} of the covering, and each $f_U$ the \emph{gate for $U$}.
\end{defn}

\begin{lem}\label{lem:gates}
Let $\bf C$ be a topological clone which has a gate covering, and let $\xi\colon \mathbf C\To \mathbf C'$ be a homomorphism to a topological clone $\mathbf C'$. If the restriction of $\xi$ to $\mathbf C^{(1)}$ is Cauchy continuous, then so is $\xi$.
\end{lem}
\begin{proof}
Let $(g^j)_{j\in\omega}$ be a Cauchy sequence in $\mathbf C$; we have to show that $(\xi(g^j))_{j\in\omega}$ is Cauchy in  $\mathbf C'$. Let $\mathcal U$ be a gate covering of $\mathbf C$. We may assume that all $g^j$ have equal arity $n\geq 1$, and that there exists $U\in\mathcal U$ containing all $g^j$. Let $\alpha^j,\beta_i^j$ be as in the definition of a gate covering. Then because the restriction of $\xi$ to $\mathbf C^{(1)}$ is Cauchy continuous, the images of the sequences $(\alpha^j)_{j\in\omega}$, $(\beta_i^j)_{j\in\omega}$ under $\xi$ are Cauchy. Hence, $(\xi(g^j))_{j\in\omega}$ is Cauchy as well, because $\xi(g^j)(x_1,\ldots,x_n)=\xi(\alpha^j)(\xi(f_U)(\xi(\beta^j_1)(x_1),\ldots,\xi(\beta^j_n)(x_n)))$ for all $j\in\omega$ and because composition is continuous.
\end{proof}

\begin{thm}\label{thm:lifthigher}
Let $\mathbf C$ be a closed subclone of $\mathbf O$ which has a gate covering, is transitive, and such that $\mathbf C^{(1)}$ has automatic homeomorphicity. Then $\mathbf C$ has automatic homeomorphicity.
\end{thm}
\begin{proof}
Let $\mathbf C'$ be another closed subclone of $\mathbf O$ and let $\xi\colon \mathbf C\To\mathbf C'$ be an isomorphism. Then $\xi$ is open by Proposition~\ref{prop:trans-open} since $\mathbf C$ is transitive and since the restriction of $\xi$ to unary functions is a homeomorphism. By Lemma~\ref{lem:gates}, $\xi$ is continuous.
\end{proof}

We will now apply Theorem~\ref{thm:lifthigher} to the \emph{Horn clone}, a well-studied (cf.~\cite{BodChenPinsker}) function clone on a countable domain which plays an important role for \emph{equality constraint satisfaction problems}~\cite{ecsps, qecsps}.

\begin{defn}\label{defn:horn}
The \emph{Horn clone} $\mathscr H$ is the smallest closed function clone on a countably infinite domain $D$ which contains all injections from finite powers of $D$ into $D$. It is easy to see that $\mathscr H$ consists of all functions of the form 
$$
f(\pi^n_{i_1}(x_1,\ldots,x_n),\ldots,\pi^n_{i_k}(x_1,\ldots,x_n)),
$$ 
where $1\leq k\leq n$, $f\colon D^k\To D$ is injective, and 
$i_1,\ldots,i_k\in\{1,\ldots,n\}$ are such that $i_1<\cdots <i_k$; in this representation of a fixed  function in $\mathscr H$, the parameters $n,k,f$, and $i_1,\ldots,i_k$ are unique. We call the functions in $\mathscr H$ \emph{essentially injective}.
\end{defn}

\begin{prop}
$\mathscr H$ has automatic homeomorphicity.
\end{prop}

\begin{proof}
We show that $\mathscr H$ has a gate covering. Let $\mathscr H^{(n)}_{i_1,\ldots,i_k}$ consist of functions of the form $f(\pi^n_{i_1}(x_1,\ldots,x_n),\ldots,\pi^n_{i_k}(x_1,\ldots,x_n))$, as in Definition~\ref{defn:horn}. The sets $\mathscr H^{(n)}_{i_1,\ldots,i_k}$ are clearly closed, and the set of all $n$-ary functions in the Horn clone can be written as a finite disjoint union of such sets. So all these sets are clopen. 

For each set $\mathscr H^{(n)}_{i_1,\ldots,i_k}$, we pick a function $$f(x_1,\ldots,x_n):=f'(\pi^n_{i_1}(x_1,\ldots,x_n),\ldots,\pi^n_{i_k}(x_1,\ldots,x_n)),$$ where $f'\colon D^k\To D$ is bijective.  It is clear that every $g\in \mathscr H^{(n)}_{i_1,\ldots,i_k}$ can be uniquely  written as
$$
g(x_1,\ldots,x_n)=\alpha(f(x_1,\ldots,x_n)),
$$
where $\alpha\in\mathscr H^{(1)}$ and $f$ is the function we picked for $\mathscr H^{(n)}_{{i_1,\ldots,i_k}}$. Now if $(g^j)_{j\in\omega}$ is a converging sequence of functions in $\mathscr H^{(n)}_{i_1,\ldots,i_k}$, and we write $g^j(x_1,\ldots,x_n)=\alpha^j(f(x_1,\ldots,x_n))$, then the $\alpha^j$ converge because the $g^j$ converge and because $f$ is onto. Hence $\mathscr H$ has a gate covering. Because $\H^{(1)}$ has automatic homeomorphicity by Theorem~\ref{thm:mainmon}, and because $\mathscr H$ is transitive, Theorem~\ref{thm:lifthigher} implies that $\mathscr H$ has automatic homeomorphicity as well.
\end{proof}

\subsubsection{Advanced gate coverings}

Often the function clone for which we want to show automatic homeomorphicity does not have a gate covering itself, but a closely related clone does. We will now refine our gate covering technique to deal with this situation.

\begin{defn}
For a clone $\frak C$ and a unary element $e$ in $\frak C$, we write $e\circ \frak C$ for the smallest subclone of $\frak C$ which contains all elements of the form $e\circ f$, where $f\in \frak C$. That is, $e\circ \frak C$ consists of elements of the form $e\circ f$ as well as the projections.
\end{defn}

\begin{lem}\label{lem:composeleft}
Let $\mathscr C, \mathscr C'$ be function clones, and let $\xi\colon \mathscr C\To \mathscr C'$ be an isomorphism. Let moreover $e\in \mathscr C$ be unary and so that both $e$ and $\xi(e)$ are injective. If the restriction of $\xi$ to $e\circ \mathscr C$ is continuous (open), then also $\xi$ is continuous (open).
\end{lem}
\begin{proof}
The mappings $\psi\colon \mathscr C\To e\circ \mathscr C$ and $\psi'\colon \mathscr C'\To \xi(e)\circ \mathscr C'$ defined by $f\mapsto e\circ f$ and $f\mapsto \xi(e)\circ f$, respectively, are continuous, injective, and open. Writing $\xi'$ for the restriction of $\xi$ to $e\circ \mathscr C$, we have $\xi=\psi'^{-1}\circ\xi'\circ \psi$, proving the lemma.
\end{proof}


\begin{lem}\label{lem:advancedGates}
Let $\mathbf C$ be a cloned subclone of $\mathbf O$ with the property that there exists $e\in\mathbf C^{(1)}$ in the closure of the invertibles of $\mathbf C^{(1)}$ such that $e\circ \mathbf C$ has a gate covering. Let $\xi\colon \mathbf C\To \mathbf C'$ be an isomorphism onto a closed subclone $\mathbf C'$ of $\mathbf O$. 
If the restriction of $\xi$ to $\mathbf C^{(1)}$ is continuous, then so is $\xi$.
\end{lem}
\begin{proof}
Since $e$ is in the closure of the invertibles of $\mathbf C^{(1)}$ we have that $\xi(e)$  is in the closure of the invertibles of $\mathbf C'^{(1)}$. Hence, $e$ and $\xi(e)$ are injective in any actions of $\mathbf C$ and $\mathbf C'$ as function clones, and Lemma~\ref{lem:composeleft} applies. Because $e\circ \mathbf C$ has a gate covering, and the restriction of $\xi$ to the unary functions in $e\circ \mathbf C$ is continuous, the restriction of $\xi$ to $e\circ \mathbf C$ is continuous by Lemma~\ref{lem:gates}. Hence, $\xi$ is continuous by Lemma~\ref{lem:composeleft}.
\end{proof}
 
We now give an example of a function clone where we use Lemma~\ref{lem:advancedGates} rather than Lemma~\ref{lem:gates} in order to prove continuity. In the following, let $N$ denote the \emph{non-edge} relation on the random graph $(V;E)$, i.e., the relation  defined on the random graph by the formula $x\neq y\wedge \neg E(x,y)$.

\begin{lem}\label{lem:gate-randomgraph}
Let $e$ be a self-embedding of the random graph $(V;E)$ whose range is co-rich. Then $e\circ (\Pol(V;E)\cap\H)$ has a gate covering.
\end{lem}
\begin{proof}
For every $n\geq 1$ we construct a gate for the set $\F_n:=e\circ(\Pol(V;E)\cap\H^{(n)}_{1,\ldots,n})$ of $n$-ary injective functions in $e\circ\Pol(V;E)$ (which is clopen subset of this clone). The proof for sets of the form $\Pol(V;E)\cap \H^{(n)}_{i_1,\ldots,i_k}$ (i.e., $n$-ary functions with a fixed set of dummy variables) is similar.

Consider the language $\tau:=\{A,B,E_A,E_B,\phi,\psi_1,\ldots,\psi_n\}$, in which $A,B$ are unary predicates, $E_A,E_B$ are binary relational symbols, $\phi$ is an $n$-ary function symbol, and $\psi_1,\ldots,\psi_n$ are unary function symbols. We define the following axioms.
\begin{enumerate}
\item The sets $A,B$ form a partition of the domain;
\item $E_A,E_B$ are the edge relation of an undirected graph without loops on $A$ and $B$, respectively;
\item $\phi$ is an injective partial function whose domain is $A^n$ and whose range is contained in $B$ such that
 $\bigwedge_{1\leq k\leq n}E_A(u_k,v_k)$ implies $E_B(\phi(u_1,\ldots,u_n),\phi(v_1,\ldots,v_n))$;
\item each $\psi_k$ is a partial function defined on the range of $\phi$ such that $\psi_k(\phi(u_1,\ldots,u_n))=u_k$ for all $(u_1,\ldots,u_n)$ in the domain of $\phi$.
\end{enumerate}
Let $\mathcal C$ be the class of all $\tau$-structures satisfying these axioms. In the following, when $\Gamma\in\mathcal C$, then we write $A^\Gamma$ and so forth for the interpretations of the symbols of $\tau$ in $\Gamma$. Every triple $(g,S,T)$ where $g\in\F_n$, $S,T\subseteq V$, and $g[S^n]\subseteq T$, gives rise to a structure $\Gamma_{(g,S,T)}$ in $\mathcal C$: $A^\Gamma:=S$ ; $B^\Gamma:=T$; the relations $E_A^\Gamma,E_B^\Gamma$ are just the appropriate restrictions of $E$; $\phi^\Gamma$ equals $g$; and the $\psi_k^\Gamma$ are defined as required by Item~(4). Conversely, every countable structure  
in $\mathcal C$ is of the form $\Gamma_{(g,S,T)}$ as above. We write $\Gamma_g$ for $\Gamma_{(g,V,V)}$.

The class of finite structures in $\mathcal C$ is a Fra\"{i}ss\'{e} class if we extend the classical definition of a Fra\"{i}ss\'{e} class as in~\cite{Hodges} to structures with partial functions; this could also be avoided by adding a distinguished element $\infty$ to the structures and working with total functions which yield $\infty$ whenever our partial functions were undefined. To see the amalgamation property of $\mathcal C$, we indicate how given finite structures $\Gamma_0, \Gamma_1, \Gamma_2\in\mathcal C$ and embeddings $s_1\colon \Gamma_0\To\Gamma_1$ and $s_2\colon \Gamma_0\To\Gamma_2$ one can build a structure $\Gamma_3$ in $\mathcal C$ and embeddings $t_1\colon \Gamma_1\To\Gamma_3$ and $t_2\colon \Gamma_2\To\Gamma_3$ such that $t_1\circ s_1=t_2\circ s_2$. We may assume that $\Gamma_0$ is a substructure of $\Gamma_1$ and $\Gamma_2$, and that $s_1,s_2$ are the identity embeddings. To begin with, $A^{\Gamma_3}$ equals $A^{\Gamma_1}\cup A^{\Gamma_2}$. For $\phi^{\Gamma_3}$, one first takes the union of $\phi^{\Gamma_1}$ and $\phi^{\Gamma_2}$, and then extends the function to $A^{\Gamma_3}$ injectively, using new elements as values;  $B^{\Gamma_3}$ will consist of  $B^{\Gamma_1}\cup B^{\Gamma_2}$ plus these new values. Finally, one extends $E_B$ on $B^{\Gamma_3}$  so that Item~(3) is satisfied, and defines the $\psi_k$ as prescribed by Item~(4). The embeddings $t_1,t_2$ are just the identity functions on $\Gamma_1$ and $\Gamma_2$, respectively. The Fra\"{i}ss\'{e} limit of $\mathcal C$ is of the form $\Gamma_{f}$ for some  $f\in \F_n$; we claim that any such $f$ is a gate for $\F_n$.

Let any $g\in \F_n$ be given. We will find self-embeddings $\alpha,\beta_1,\ldots,\beta_n$ of $(V;E)$ such that $g(x_1,\ldots,x_n)=\alpha(f(\beta_1(x_1),\ldots,\beta_n(x_n))$. 
It will be convenient to construct $\beta:=(\beta_1,\ldots,\beta_n)$ as a function from $V^n$ into itself which respects equality of coordinates.

$$
\begin{CD} V^n	@>g>> V\\ @V \beta VV @AA \alpha A\\ V^n	@>f>> V \end{CD}
$$

We will construct $\beta$ and $\alpha$ simultaneously by a back-and-forth argument using the universality and homogeneity of $\Gamma_{f}$: each $\beta_k$ will be an embedding from the structure $(A^{\Gamma_g};E_A^{\Gamma_g})=(V;E)$ into $(A^{\Gamma_f};E_A^{\Gamma_f})=(V;E)$, and $\alpha$ will be an embedding from $(B^{\Gamma_f};E_B^{\Gamma_f})=(V;E)$ into $(B^{\Gamma_g};E_B^{\Gamma_g})=(V;E)$; moreover, we will have $\phi^{\Gamma_g}=\alpha\circ \phi^{\Gamma_f}\circ \beta$. Suppose that $\beta$ is already defined on $\{p^1,\ldots,p^r\}\subseteq  V^n$, and $\alpha$ is already defined on $\{v^1,\ldots,v^w\}\subseteq V$. 
We will show that we can extend the domains of $\beta$ and $\alpha$ whilst staying consistent (i.e., whilst maintaining the above conditions).

Consider first the case where we wish to extend the domain of $\beta$ to some $p\in V^n$. Then let $p'\in V^n$ be so that 
\begin{enumerate}
\item the type of $p'$ over $(\beta(p^1),\ldots,\beta(p^r),v^1,\ldots,v^w)$ in $\Gamma_f$ equals the type of $p$\\ over $(p^1,\ldots,p^r,\alpha(v^1),\ldots,\alpha(v^w))$ in $\Gamma_g$ (in particular, $\beta$ remains an embedding).
\end{enumerate}
We set $\beta({p}):=p'$ and $\alpha(f(p')):=g({p})$. Easy!

Next consider the case where we wish to extend the domain of $\alpha$ to some $v\in V$. If $v$ is contained in the range of $f$, and for the unique element $p'\in V^n$ with $f(p')=v$ we have that there exists $p\in V^n$ such that the condition above is valid, then we can simply extend $\beta$ to $p$ as above, hence thereby extending $\alpha$ to $v$. However, $v$ may either not be contained in the range of $f$, or $p$ might not exist, because $\Gamma_g$ is not necessarily universal or homogeneous; in this case, we have to use the fact that we can send $v$ to an element outside the range of $g$. It is here that we use the co-richness of the image of $g$. In this situation, let $v'\in V$ be so that
\begin{enumerate}
\item[(2)] $v'$ is not contained in the range of $g$;
\item[(3)] the type of $v'$ over $(\alpha(v^1),\ldots,\alpha(v^w))$ in $(V;E)$ equals the type of $v$ over $(v^1,\ldots,v^w)$ in $(V;E)$ (i.e., $\alpha$ remains an embedding).
\end{enumerate}
Condition~(2) ensures that we do not get stuck in the future when we wish to extend $\beta$ to some $p\in V^n$ with $g({p})=v$. Now set $\alpha(v):=v'$.

The fact that these two steps are always possible proves the gateness of $f$ for $\F_n$, since we can extend the embeddings $\alpha$ and $\beta$ in an alternating fashion to obtain total functions.

Now let $(g^j)_{j\in\omega}$ be a converging sequence of functions in $\F_n$. By the above, we can write every $g^j$ as
$$
g^j(x_1,\ldots,x_n)=\alpha^j(f(\beta^j_1(x_1),\ldots,\beta^j_n(x_n))),
$$
for self-embeddings $\alpha^j,\beta^j_1,\ldots,\beta^j_n$ of $(V;E)$. Write $\beta^j:=(\beta^j_1,\ldots,\beta^j_n)$. Now for every finite set $F\subseteq V^n$ there exists $d_F\in\omega$ such that all $g^j$ with $j\geq d_F$ agree on $F$. We can start the construction of $\beta^j$ assuming that $\beta^j$ agrees with $\beta^{d_F}$ on $F$, since our construction only depends on what $g^j$ does on $F$, and not on what it does outside. Hence, we can assure that the $\beta^j$ agree on $F$ from $d_F$ on, and thus assure that $(\beta^j)_{j\in\omega}$ converges. In the same way, for every finite set $S\subseteq V$, we can start our construction of $\alpha^j$ and of $\beta^j$ with any given values for $\alpha^j$ on $S$. We can thus also assure that $(\alpha^j)_{j\in\omega}$ converges.
\end{proof}

We immediately obtain the following conditional corollary stating that continuity lifts from the unary functions in $\Pol(V;E)\cap\H$ to the entire clone for isomorphisms from this clone. We will shortly see, via a similar technique, that requiring continuity on the unary functions is not even necessary (Lemma~\ref{lem:endrandom}).

\begin{cor}\label{cor:ah-randomgraph}
Consider the function clone $\Pol(V;E)\cap\H$ of essentially injective functions in $\Pol(V;E)$, and let $\xi$ be any isomorphism from this clone onto another closed function clone on a countable set. If the restriction of $\xi$ to the unary functions is continuous, then $\xi$ is continuous.
\end{cor}
\begin{proof}
This is a direct consequence of Lemmas~\ref{lem:gate-randomgraph} and~\ref{lem:advancedGates}.
\end{proof}

\begin{thm}\label{thm:lifthigher-improved}
Let $\mathbf C$ be a closed subclone of $\mathbf O$ which is transitive, and such that $\mathbf C^{(1)}$ has automatic homeomorphicity. If there exists $e\in\mathbf C^{(1)}$ in the closure of the invertibles of $\mathbf C^{(1)}$ such that $e\circ \mathbf C$ has a gate covering, then $\mathbf C$ has automatic homeomorphicity.
\end{thm}
\begin{proof}
Let $\xi\colon\mathbf C \To \mathbf C'$ be an isomorphism onto a closed subclone $\mathbf C'$ of $\mathbf O$. Then $\xi$ is continuous by Lemma~\ref{lem:advancedGates} and open by Proposition~\ref{prop:trans-open}.
\end{proof}

There exist 17 closed function clones containing the automorphism group of the random graph such that a relational structure with a first-order definition in $(V;E)$ has a tractable CSP if and only if its polymorphism clone contains one of the function clones of that list~\cite{BodPin-Schaefer}. The clones of that list are called  \emph{minimal tractable polymorphism clones} over the random graph. With one exception, they all consist of essentially injective functions with a certain predescribed behavior with respect to edges and non-edges. The method of Lemma~\ref{lem:gate-randomgraph} yields the following.

\begin{cor}\label{cor:minimaltractable}
All minimal tractable polymorphism clones over the random graph have automatic homeomorphicity.
\end{cor}
\begin{proof}
The exception among the 17 function clones mentioned above is the function clone corresponding to the transformation monoid consisting of $\overline{\Aut(V;E)}$ plus all constant functions. If $\xi$ is an isomorphism from this clone onto another closed subclone of $\mathscr O_V$, then the restriction of $\xi$ is to $\overline{\Aut(V;E)}$ is continuous and open by Corollary~\ref{cor:mon}; since constant functions are sent to constant functions, it is clear that $\xi$ is continuous and open as well.

The other 16 function clones all have $\overline{\Aut(V;E)}$ as their unary part and consist of essentially injective functions, similarly to $\Pol(V;E)\cap\H$. The only difference with the latter clone is that there are stronger requirements on the functions in the clone than 
preservation of $E$ (for example, in all cases they also preserve $N$). Adjusting Item~(3) in the axioms of the Fra\"{i}ss\'{e} class accordingly, the very same proof as in Lemma~\ref{lem:gate-randomgraph} goes through to show automatic homeomorphicity via Theorem~\ref{thm:lifthigher-improved}.
\end{proof}

\subsubsection{Gates for endomorphisms and decomposing functions: $\Pol(V;E)$}\label{sect:gatesforunary}

One can also use the gate technique in order to show continuity of, say, homomorphisms from the endomorphism monoid of a structure given continuity of their restrictions to the self-embedding monoid. Consider, for example, the random graph $(V;E)$; then its self-embedding monoid $\overline{\Aut(V;E)}$ has automatic homeomorphicity by Corollary~\ref{cor:mon}. Using the same method as in Lemma~\ref{lem:gate-randomgraph}, we can construct an endomorphism $f$ of $(V;E)$ such that for all converging sequences $(g^j)_{j\in\omega}$ of functions in $e\circ \End(V;E)$, where $e$ is a fixed self-embedding of $(V;E)$ with co-rich range, we have self-embeddings $\alpha^j,\beta^j$ of $(V;E)$ such that
\begin{itemize}
\item $g^j=\alpha^j\circ f\circ \beta^j$ for all $j\in\omega$;
\item $(\alpha^j)_{j\in\omega}$ and $(\beta^j)_{j\in\omega}$ converge.
\end{itemize}

\begin{lem}\label{lem:endrandom}
Let $\xi$ be an isomorphism from $\End(V;E)$, or from $\End(V;E)\cap\H$, onto another closed submonoid of $\mathscr O^{(1)}_V$. Then $\xi$ is continuous.
\end{lem}
\begin{proof}
When constructing the gate $f$ as in Lemma~\ref{lem:gate-randomgraph}, we of course set $n=1$, and do not require $\phi$ to be injective in the case of  $\End(V;E)$ (since we do not have to maintain injectivity in the amalgamation, we can also drop Item~(4)).
The rest of the proof is identical with the proof of Lemma~\ref{lem:gate-randomgraph}.
\end{proof}

\begin{lem}\label{lem:randomgraph:continuous}
Let $\xi$ be an isomorphism from $\Pol(V;E)$ onto another closed subclone of $\mathscr O^{}_V$. Then $\xi$ is continuous.
\end{lem}
\begin{proof}
We know that the restrictions of $\xi$ to $\End(V;E)$ and to $\Pol(V;E)\cap\H$ are continuous by Corollary~\ref{cor:ah-randomgraph} and Lemma~\ref{lem:endrandom}. Let $e$ be a self-embedding of the random graph $(V;E)$ whose range is co-rich. We show that the restriction of $\xi$ to $e\circ \Pol(V;E)$ is continuous; the theorem then follows from Lemma~\ref{lem:composeleft}.

We first show that if $g$ is any $n$-ary function in $e\circ \Pol(V;E)$ for some fixed $n\geq 1$, then $g$ can be written as $h\circ f$, where $f\in\Pol(V;E)$ is an $n$-ary injective function, and $h\in \End(V;E)$. 

To this end, suppose that values under $h$ and $f$ have already been defined for $v^1,\ldots,v^w\in V$ and for  $p^1,\ldots,p^r\in V^n$, respectively. 

If we wish to define $f$ on some $p\in V^n$, then pick $f({p})\in V\setminus\{v^1,\ldots,v^w\}$ such that $E(f({p}),v^i)$ iff $E(g({p}),h(v^i))$, for all $1\leq i\leq w$; subsequently, set $h(f({p})):=g({p})$. Clearly, that way $h$ remains a partial self-embedding of $(V;E)$. To see that the extended function $f$ still preserves $E$, suppose that $p$ and $p^j$ are adjacent in all components. Then $E(g({p}),g(p^j))$, or written differently, $E(g({p}),h\circ f(p^j))$
 holds. Writing $v^i=f(p^j)$ we see that our choice of $f({p})$ implies $E(f({p}),v^i)$, and so $E(f({p}),f({p^j}))$.

 If we wish to extend the domain of $h$ to some $v\in V$, then simply pick $h(v)$ outside the range of $g$ so that it remains an embedding.

Given a sequence $(g^j)_{j\in\omega}$ of functions in $e\circ \Pol(V;E)$, we can write each $g^j$ as $h^j\circ f^j$ as above. Moreover, we can pick the $h^j$ and $f^j$ so that they converge: for from the construction of those functions we see that if they are (consistently) predefined on finite domains, then one can always extend them and complete the construction. Since $(h^j)_{j\in\omega}$ and $(f^j)_{j\in\omega}$ converge, so do $(\xi(h^j))_{j\in\omega}$ and $(\xi(f^j))_{j\in\omega}$, because the restrictions of $\xi$ to $\End(V;E)$ and to $\Pol(V;E)\cap\H$ are continuous by Corollary~\ref{cor:ah-randomgraph} and Lemma~\ref{lem:endrandom}. Hence, $(\xi(g^j))_{j\in\omega}=(\xi(h^j)\circ \xi(f^j))_{j\in\omega}$ converges to $\xi(h)\circ\xi(f)=\xi(h\circ f)=\xi(g)$, where $h,f,g$ are the limits of the sequences $(h^j)_{j\in\omega}$, $(f^j)_{j\in\omega}$, and $(g^j)_{j\in\omega}$, respectively.
\end{proof}

\subsection{Topological Birkhoff and openness}\label{subsect:topobirkopenness}

When $\mathcal C$ is a class of algebras with common signature $\tau$, then $\PPPfin(\mathcal C)$ denotes the class of all finite products of algebras from $\mathcal C$. The following is a variant of Theorem~\ref{thm:birkhoff} for certain function clones.

\begin{thm}[Bodirsky and Pinsker~\cite{Topo-Birk}]\label{thm:topobirkhoff}
Let $\mathscr C, \mathscr D$ be closed function clones acting on countable sets $C$, $D$ respectively such that $\mathscr C$ is oligomorphic and such that the algebra $(D;\mathscr D)$ is finitely generated. Then there exists a surjective continuous homomorphism $\xi\colon \mathscr C\To \mathscr D$ if and only if $(D;\mathscr D)\in \HSPfin(C;\mathscr C)$ for some indexing of those algebras.
\end{thm}

\begin{lem}\label{lem:randomgraph:open}
Let $\xi\colon \Pol(V;E)\To\mathscr D$ be an isomorphism onto a closed subclone $\mathscr D$ of $\mathscr O_V$. Then $\xi$ is open.
\end{lem}
\begin{proof}
Consider any finitely generated subalgebra $B$ of $(V;\mathscr D)$ with at least two elements.
By Lemma~\ref{lem:randomgraph:continuous} we know that $\xi$ is continuous, and hence so is the mapping $\xi'$ which sends every function $f\in\Pol(V;E)$ to $\xi(f)\rest_B$. Hence, by Theorem~\ref{thm:topobirkhoff}, the algebra $(B;\mathscr D)$ is a homomorphic image of a subalgebra of the algebra $(V; \Pol(V;E))^n$, for some $n\geq 1$. We write $S$ for the universe of that subalgebra, and  $\sim$ for the congruence relation on $S$ such that $(S; \Pol(V;E))/_\sim$ is isomorphic to $(B;\mathscr D)$. We will show that $\xi'$ is open and injective. Then $\xi$ is open as well, because $\xi=r^{-1}\circ\xi'$ for the continuous mapping $r$ which restricts every function in $\xi[\Pol(V;E)]$ to $B$.

Let $Q^E$ ($Q^=$) consist of all pairs $(i,j)$ such that $1\leq i,j\leq n$ and such that $E(a_i,a_j)$ ($a_i=a_j$) for all $a\in S$. Write $P$ for those $(i,j)$ with $1\leq i,j\leq n$ which are contained in neither of the two sets. Then there exists a tuple $t\in S$ such that $N(t_i,t_j)$ for all $(i,j)\in P$: it is enough to pick for every pair $(i,j)\in P$ tuples $t_{i,j}^E$ and $t_{i,j}^=$ witnessing that $(i,j)\nin Q^E$ and $(i,j)\nin Q^=$, and then apply a function $f\in\Pol(V;E)$ to all those tuples so that the result is a tuple $t$ as desired. Now it is clear that for any two tuples $a,c$ in $S$, there exists a binary function $f\in\Pol(V;E)$ such that $f(a,t)=c$: setting the values on $(a,c)$ does not violate $E$, and $f$ can be extended to a total function by the universality and homogeneity of the random graph. For $a,b\in S$, write $I(a,b)$ for the set of all $1\leq i\leq n$ such that $a_i= b_i$. Then if $a,b,c,d\in S$ are so that $I(a,b)\subseteq I(c,d)$, then $f(a,t)=c$ and $f(b,t)=d$ for an appropriate function $f\in\pol(V;E)$, by the same argument as above. Hence, in this situation $a\sim b$ implies $c\sim d$. In particular, whether or not $a\sim b$ holds for given $a,b\in S$ only depends on $I(a,b)$. Let $W$ be the set of all subsets of $\{1,\ldots,n\}$ of the form $I(a,b)$, where $a,b\in S$ and  $a\sim b$. Then by the above, $W$ is upward closed and for $a,b\in S$ we have $a\sim b$ if and only if $I(a,b)\in W$. Moreover, 
$W$ is closed under intersections: when $I,J\in W$, then there exist tuples $a,b,c\in S$ such that $I(a,b)=I$, $I(b,c)=J$, and $I(a,c)=I\cap J$. Since $a\sim b$ and $b\sim c$ imply $a\sim c$ we infer $I\cap J\in W$. We cannot have $\emptyset\in W$, for otherwise $W$ would contain all subsets of $\{1,\ldots,n\}$, and hence $\sim$ would identify all tuples of $S$; but this contradicts our assumption that $B$ has more than one element. Consequently, the intersection of all sets in $W$ is non-empty. Picking any $1\leq i\leq n$ in this intersection, we have that $a_i\neq b_i$ implies that $a\sim b$ does not hold, for all $a,b\in S$. Moreover, because $\Aut(V;E)$ acts transitively on $V$, the projection of the set $S$ onto its $i$-th coordinate equals $V$.

It follows immediately that $\xi'$ is injective. Let $U$ be a basic open subset of $\Pol(V;E)$, i.e., $U$ consists of all $k$-ary $f\in\Pol(V;E)$ which satisfy $f(a_1,\ldots,a_k)=a_0$, for fixed $a_0,\ldots,a_k\in V$. Then picking any $\sim$-classes $A_0,\ldots,A_k$ such that $t\in A_j$ implies $t_i=a_i$ for all $1\leq j\leq k$, the set $U$ can be described via the action of $\Pol(V;E)$ on $B$ as follows: it consists of those $k$-ary $f\in\Pol(V;E)$ which satisfy $f(A_1,\ldots,A_k)=A_0$. Hence, $\xi[U]$ is clopen.
\end{proof}

\begin{thm}
$\Pol(V;E)$, the polymorphism clone of the random graph, has automatic homeomorphicity.
\end{thm}
\begin{proof}
This follows from Lemmas~\ref{lem:randomgraph:continuous} and~\ref{lem:randomgraph:open}.
\end{proof}

\section{Open Problems}\label{sect:open}

It is known that the closed subgroups of ${\bf S}$ 
are precisely those topological groups that are Polish
and have a left-invariant ultrametric~\cite{BeckerKechris}. 

\begin{quest}\label{prob:topo-monoids}
Give a characterization of those topological monoids that appear as closed submonoids of ${\bf O}^{(1)}$.
\end{quest}

\ignore{
\begin{quest}\label{quest:monoid-topology}
Let $\bf M$ be a topological monoid.
 Are the following equivalent?
\begin{itemize}
\item $\bf M$ is a closed submonoid of ${\bf O}^{(1)}$.
\item The topology of ${\bf M}$ is Polish and $\bf M$ has a left nonexpansive complete ultrametric.
\end{itemize}
\end{quest}

We can give a positive answer to Question~\ref{quest:monoid-topology} in the case where
the invertibles are dense in $\bf M$
(for the interesting direction, we apply the result for groups to the group $\bf G$ of invertibles in $\bf M$.
We obtain a structure $\Gamma$ whose automorphism group is $\bf G$; this structure
$\Gamma$ has $\bf M$ as endomorphism monoid).

A compatible metric $d$ on 
a topological clone ${\bf C}$ 
is called \emph{non-expansive} if 
for all $f \in C^{(k)}$ and $g_1,g_1',\dots,g_k,g_k' \in C^{(k)}$,
we have that 
$$d(\comp^k_k(f,g_1,\dots,g_k),\comp^k_k(f,g'_1,\dots,g'_k)) \leq \max(d(g_1,g_1'),\dots,d(g_k,g_k')) \; .$$
Note that the complete metric on function clones defined above is a left nonexpansive ultrametric.
}

\begin{quest}
Give a characterization of those topological clones that appear as closed subclones of ${\bf O}$.
\end{quest}

In connection with Theorem~\ref{thm:no-auto-homeo} and the known counterexample for closed oligomorphic groups~\cite{EvansHewitt}, we ask the following.

\begin{quest}
Is there a closed oligomorphic subclone of $\bf O$ which does not have reconstruction?
\end{quest}







The following is an example of a relatively simple function clone where our techniques fail.

\begin{quest}
Does $\Pol(\mathbb Q;<)$ have automatic homeomorphicity?
\end{quest}

{\bf Recent Progress.} All four questions that have
been asked in the submitted version have been solved by the time of the preparation of the revised version of the article. For Questions~1 and~2,
this appeared in recent work with F. Martin Schneider\cite{BodirskySchneiderTopological}, including a topological characterisation of the topological clones that appear as \emph{oligomorphic} closed subclones of ${\bf O}$. A positive answer to Question~2
is the main result of recent work with David Evans and Michael Kompatscher~\cite{BodirskyEvansKompatscherPinsker}. 
Finally, a positive answer to the variant
of Question~4 where $(\mathbb Q;<)$ is replaced by $({\mathbb Q};\leq)$ 
has been announced by Mike Behrisch, 
John Truss, and Edith Vargas-Garc\'{i}a~\cite{BehrischTrussVargas}. 
For more results of automatic homeomorphicity 
of clone homomorphisms we refer to work of Christian Pech and Maja Pech~\cite{PechPechHomeo}. 
Finally, a positive answer to Question~4 follows from more general results of Robert Barham~\cite{BarhamLocallyMoving}. 




\section*{Acknowledgements}
We are grateful to John Truss, Edith Vargas-Garc\'{i}a,
and Christian Pech for their comments on earlier versions of the paper. The first author thanks Dietrich Kuske for pointing him in 2003 to the work of Rubin.

\bibliographystyle{alpha}
\bibliography{7_reconstruction.bib}

\end{document}